\documentclass[12pt]{elsarticle}
\usepackage{graphicx}
\usepackage{numcompress}
\usepackage{amssymb,amsmath,amsthm}
\usepackage{epstopdf}
\usepackage{pstricks}

\theoremstyle{plain}
\newtheorem{thm}{Theorem}  
\newtheorem{cor}[thm]{Corollary}
\newtheorem{lem}[thm]{Lemma}

\newtheorem*{cor*}{Corollary}
\newtheorem*{con*}{Conjecture}

\newtheorem*{thm*}{Theorem}
\newtheorem{lem*}{Lemma}
\newtheorem{prop}[thm]{Proposition}
\newtheorem*{prop*}{Proposition}
\newtheorem{defn}{Definition} 
\newtheorem*{defn*}{Definition}  
\newtheorem*{obsrv*}{Observations} 
\newtheorem{rem}{Remark}   
 \theoremstyle{definition}
\newtheorem*{open*}{Open Problem}

\newcommand\bea{\begin{eqnarray}}
\newcommand\eea{\end{eqnarray}}
\newcommand\bi{\begin{itemize}}
\newcommand\ei{\end{itemize}}
\newcommand\ben{\begin{enumerate}}
\newcommand\een{\end{enumerate}}
\newcommand{\comments}[1]{}

\numberwithin{equation}{section}

\begin{document}
\begin{Large}
\begin{center}
\bf{Microlocal analysis of scattering data for nested conormal potentials}
\end{center}
\begin{center}
Suresh Eswarathasan
\end{center}
\end{Large}

\begin{footnotesize}

\noindent Abstract: Working in the time domain, we show that the location of the singularities and the principal symbol of a potential that is conormal to nested submanifolds $S_2 \subset S_1 \subset \mathbb{R}^n$, for $n \geq 3$, can be recovered from the backscattering as well as from the restriction of the far-field pattern to more general determined sets of scattering data.  This extends the work of Greenleaf and Uhlmann where the potentials considered are conormal to a single submanifold $S \subset \mathbb{R}^n$.  We utilize the microlocal analysis of the wave operator $\square=\partial_t^2 - \triangle_x$ and multiplication by a nested conormal distribution in order to study their action on spaces of conormal-like distributions.
\end{footnotesize}

\section{Introduction} \label{intro}

Consider the potential scattering problem for the wave equation:
\begin{eqnarray} \label{potprob}
\nonumber &&(\partial_t^2 - \triangle + q)u=0 \text{ in } \mathbb{R}^n \times \mathbb{R}, \\ 
&& u=\delta(t - x \cdot \omega) \text{ for } t << - \rho, 
\end{eqnarray}
where $q$ is a compactly supported time-independent potential and $\omega \in S^{n-1}$ is varying.  Here $\rho$ is any value such that $supp(q) \subset \{|x| \leq \rho\}$.  One obtains a scattering kernel for a fixed $q$, say in $\in C^{\infty}_0(\mathbb{R}^n)$, from the Friedlander radiation condition \cite{MelUhlnotes},
\begin{equation}
\alpha_q(t, \theta, \omega) = \lim_{r \rightarrow \infty} r^{\frac{n-1}{2}}\partial_t u(t + r, r\theta, \omega),
\end{equation}
where $(r,\theta)$ are polar coordinates.  The scattering map $\Phi$ which sends $q$ to $\alpha_q$ is nonlinear and overdetermined and there has been much interest in the inverse problem of determining $q$ from $\alpha_q$.  Since $\alpha_q$ is overdetermined, it is naturally also of interest to try to reconstruct $q$ from the restriction of $\alpha_q$ to various submanifolds of $\mathbb{R} \times S^{n-1} \times S^{n-1}$.  We write the scattering amplitude for the analogous stationary potential scattering problem as $a_q(\lambda, \theta, \omega)$. 

\subsection{Statement of the problem and main result} \label{spmr}

The class of $q$'s to be considered in (\ref{potprob}) are those that have singularities conormal \cite{Ho71} to a nested pair of submanifolds $S_2 \subset S_1$ of $\mathbb{R}^n$, denoted by $(S_1, S_2)$, of arbitrary codimension.  The inverse problem we solve consists of determining these submanifolds and the principal symbol of $q$, which is enough to determine the singularities of $q$, from the leading singularities of the backscattering $\alpha_{|\mathbb{B}}$, where $\mathbb{B} = \{\theta = - \omega\} \subset \mathbb{R} \times S^{n-1} \times S^{n-1}$, which is a distribution on $\mathbb{R} \times S^{n-1}$.  We treat similar determined sets of scattering data.  

It is shown that $\alpha$, away from $\omega$'s that are tangent to either of the submanifolds, is the sum of a paired Lagrangian distribution associated to two cleanly intersecting reflected Lagrangians, two reflected Lagrangian distributions, and a single peak Lagrangian distribution, modulo Sobolev errors.  Although the strongest singularity lies on the peak Lagrangian, as is well known in the physics literature, we show that it is the restriction of the reflected Lagrangians and their points of intersection to various submanifolds of scattering data in $\mathbb{R} \times S^{n-1} \times S^{n-1}$ that determine the singularities of $q$.  The precise theorem is the following:

\begin{thm}\label{totalresult}
Let $S_2 \subset S_1 \subset \mathbb{R}^n$ be smooth nested submanifolds of codimension $d_1 + d_2$ and $d_1$, respectively.  Assume that $q$ is compactly supported and is conormal to the nested pair $(S_1, S_2)$ of orders $M_1$ and $M_2$. Furthermore, suppose that
\begin{eqnarray*}
&& M_2 > -d_2 \text{ and } M_1 < -d_1 - \frac{d_2}{2} + 1 \text{ or } \\
&& M_2 \leq -d_2 \text{ and } M_1 < -d_1 + 1, \text{ with } \\
M_1 + \frac{M_2}{2} &<& \inf\{-\frac{n-2}{n}(d_1+d_2), -d_1 - d_2 + 1\} \text{ if }  n \geq 5,
\end{eqnarray*}
and
\begin{eqnarray*}
M_1 + \frac{M_2}{2} &<& \inf\{-\frac{d_1+d_2}{2}, -d_1 - d_2 + 1\}\text{ if } n=3 \text{ and }4.
\end{eqnarray*}
Then $S_1, S_2$, and the principal symbol of $q$ are determined by the singularities of $\alpha$ restricted to the backscattering surface $\{\theta = - \omega\} \subset \mathbb{R} \times S^{n-1} \times S^{n-1}$.
\end{thm}
\noindent In fact, we prove a stronger result that shows that for a given submanifold of scattering data with certain geometrical properties, the inverse problem can still be solved.  See Section \ref{sip} for a more detailed statement.

The orders $-\frac{n-2}{n}(d_1+d_2)$ and $-\frac{d_1+d_2}{2}$ in Theorem \ref{totalresult} are present so that we can follow the scattering theory of Lax and Phillips for short range potentials in certain $L^P(\mathbb{R}^n)$ spaces; see \cite{Phil82}.  The restriction of $ M_1 + \frac{M_2}{2} < -d_1 - d_2 + 1$ is needed in order to have a series describing singularities, which appears in a later section, stabilize in an appropriate sense.  Moreover, the potentials under consideration are allowed to blow up.  Specifically, for $M_2 < -d_2$, $q$ blows up on $S_1$ and is continuous on $S_2$, while for $0< -d_1 - \frac{d_2}{2} + 1 - M_1 < \varepsilon$ and $0< M_2 - d_2 < \varepsilon$, $q$ blows up on $S_2$ and is continuous on $S_1 \backslash S_2$.  Hence, there is no size restriction on $q$, in contrast to, e.g., \cite{Ruiz01}.

Potentials that are conormal to a single submanifold were dealt with by Greenleaf and Uhlmann \cite{GU93};  this paper closely follows the time-dependent approach taken in \cite{GU93} and generalizes the results to nested $q$.  Using Proposition 3.7 in \cite{GuiUhl81} involving the intersection of classes of paired Lagrangians over all orders $M_2$, it follows that Theorem \ref{totalresult} covers the main result of \cite{GU93}.

A significant difference between the work in this paper and that of \cite{GU93} is the new, more complicated geometry that arises when using an approximation method, the understanding of multiplication by $q$ on Sobolev spaces and other classes of distributions, and the appearance of distributions that are associated to cleanly intersecting triples and quadruples of Lagrangians.    

We conclude this section by noting that even for an arbitrary Lagrangian distribution $u$, calculating the blowup rates that assist in finding which $L^p$ space $u$ belongs to is difficult without some additional assumptions on the Lagrangian.  Hence, from the viewpoint of the Lax-Phillips scattering theory, assuming that $u$ is in some conormal category is a reasonable restriction.

\subsection{Previous results}

\noindent 1) \underline{Fixed angle scattering}:  One wants to determine information about the potential $q$ from the scattering amplitude restricted to a fixed incident angle $\theta_0 \in S^{n-1}$.   Stefanov \cite{Stef92} proves uniqueness of the potentials under a smallness assumption and Ruiz \cite{Ruiz01} shows that the Born approximation determines a ``close" approximation of $q \in H^s(\mathbb{R}^n)$ for $n=2$ and 3.

\noindent 2) \underline{Fixed energy scattering}: Here we set $\lambda = \lambda_0$ in the scattering amplitude.  Uniqueness and reconstruction results are obtained by Nachman \cite{Nach92}, Novikov \cite{Nov88}, and Ramm \cite{Ramm98} for dimensions $n \geq 3$.  Sun and Uhlmann \cite{SunUhl93} show the recovery of singularities for $n=2$.  In this setting, the main technique seems to rely on the complex exponential method used in \cite{SylUhl87} and the $\bar{\partial}$-method utilized by Nachman \cite{Nach96}.

\noindent 3) \underline{Backscattering}: As mentioned earlier in Section \ref{intro}, we set $\theta = - \omega$ in the scattering amplitude.  Uniqueness under a smallness assumption is obtained by Lagergen \cite{Lag01} for $n=3$ and recovery of singularities for $n=2$ by Ola et. al \cite{OlaPaiSer01}.  Generic uniqueness is proven in Eskin and Ralston \cite{EskRals89}.

\subsection{Outline of the paper}
This work deals with the inverse scattering problem in the time domain and is done in order to utilize a more geometrical approach through the microlocal analysis of the wave equation and multiplication by conormal-like distributions, which allows precise tracking of singularities.   

Section \ref{pre} gives the necessary background from microlocal analysis with an emphasis on Lagrangian distributions, Fourier integral operators, and paired Lagrangians distributions.  Here, we introduce the space of nested conormal distributions for which our potential $q$ will be an element.

In Section \ref{ntc}, we develop a new kind of Lagrangian-type distribution that will be seen in Sections \ref{mnc}, \ref{pwo}, and \ref{bose}.  This work is a fairly straightforward generalization of the theory of nested conormal distributions.

Sections \ref{lscr} contains a number of geometrical computations and lemmas necessary in solving our inverse problem.

The main calculations in solving the forward problem begin in Sections \ref{mnc} and \ref{pwo} where we try to understand the mapping properties of two operators, namely multiplication by $q$ and the parametrix to the wave equation.  Section \ref{bose} finishes the description of the leading singularities of the solution $u$ to (\ref{potprob}) using the Born series, relying heavily on the results from the previous sections.

The inverse problem is solved in Section \ref{sip}, where we describe the scattering kernel $\alpha_q$ using the Lax-Phillips scattering theory \cite{LP89} and the Born series from Section \ref{bose}.  The paper concludes by showing that the various restrictions of the scattering kernel continue to determine the singularities of $q$.

The paper is a revision of the author's Ph.D. thesis at the University of Rochester.

\section{Preliminaries} \label{pre}

Throughout this paper, unless otherwise specified, $X$ will represent a smooth manifold of dimension $n$, $(T^*(X) \backslash 0, \omega_{T^*(X)})$ will be the cotangent bundle of $X$ with 0-section deleted, equipped with the canonical symplectic 2-form $\omega_{T^*(X)} = \Sigma_i d \xi_i \wedge d x_i$.  See \cite{Duis11} for an overview of symplectic geometry.  The relation $a\lesssim b$ denotes $a \leq C b$ for some constant $C>0$, which may depend on some parameters but not those of interest.

\subsection{Fourier integral distributions}

Let $S \subset X$ be a smooth submanifold of codimension $k$.  Suppose $S = \{x \in X : h_j(x)=0, 1 \leq j \leq k\}$ is a local representation of $S$ with $\{\nabla h_j: 1 \leq j \leq k\}$ linearly independent on $S$.
\begin{defn}
We say that a distribution $q$ is \underline{conormal to $S$ of order $\mu$}, $\mu  \in \mathbb{R}$, if
\begin{equation}
q(x) = \int_{\mathbb{R}^k} e^{i \sum_j h_j(x) \theta_j} a(x;\theta) d \theta,
\end{equation}
with $a(x; \theta) \in S^{\mu}(X \times (\mathbb{R}^k \backslash 0))$ 
and the space of all these is denoted by $I^{\mu}(S)$.   We call  $a$ the \underline{symbol} of $q$ where $|\partial_x^{\gamma}\partial_{\theta}^{\alpha} a(x; \theta)| \lesssim  \langle \theta \rangle ^{p - |\alpha|}$ for $(x,y) \in K$ where the constant associated to $\lesssim$ depends only on $K, \alpha,$ and $\gamma$.  
\end{defn}
\noindent Here, we use the standard notation $\langle \theta  \rangle = (1 + |\theta|^2)^{\frac{1}{2}}$, the Japanese bracket of $\theta$. 

If $-k < \mu < 0$ then $q$ has a specific blowup rate: $|q(x)| \leq C \cdot dist(x,S)^{-k - \mu}$; see Section 6, Section 4 of \cite{Ste93} for the calculation of this estimate. Examples of such $q$'s are the surface measure $dS$ along $S$ as well as a function that has a Heaviside singularity across $S$.  These distributions lie in $I^{0}(S)$ and $I^{-1}(S)$, respectively.

For the conormal distributions $q$ associated to a submanifold $S$, $WF(q) \subset N^*(S)$ where $N^*(S)$ is the conormal bundle of the $S$.  These, and many other kinds of geometric distributions, fall under the ubiquitous category of \underline{Lagrangian distributions}; an important part of microlocal analysis is the study of such distributions.  For the basic theory of conic and Lagrangian submanifolds of $T^*(X)$ and nondegenerate phase functions parametrizing Lagrangian submanifolds, see \cite{Ho71}.
\begin{defn} 
Let $\Lambda \subset T^*(X) \backslash 0$ be conic Lagrangian submanifold.  The H\"{o}rmander space $I^m(\Lambda)$ of \underline{Lagrangian distributions associated to $\Lambda$} consists of locally finite sums of distributions of the form 
\begin{equation}
u(x) = \int_{\mathbb{R}^k} e^{i\varphi(x;\theta)} a(x; \theta) d \theta,
\end{equation} 
where $\varphi$ is a nondegenerate phase function parametrizing $\Lambda$, $a \in S^{m+\frac{n}{4}-\frac{k}{2}}$, and $WF(u) \subset \Lambda$.
\end{defn}
As mentioned previously, distributions conormal to a submanifold $S$ fall into the class of Lagrangian distributions with $\Lambda = N^*(S)$.  By definition, $I^{\mu}(S) = I^{\mu-\frac{n}{4}+\frac{k}{2}}(N^*(S))$.  The work of Melrose \cite{Melnotes89} gives an alternate characterization of the conormal distributions.
\begin{thm*} \label {iterreg1}
Let $S \subset X$ be a smooth submanifold.  The space of conormal distributions on $X$ with respect to $S$ of order $m$, $I^{m}(X;S)$, can be characterized as the set of all distributions $u \in \mathcal{D}'(X)$ such that 
\begin{equation*}
V_1...V_ku \in H^{-m-\frac{n}{4}, \infty}(X),
\end{equation*}
where the $V_j$'s are $C^{\infty}$ vector fields on $X$ which are tangent to $S$, and $H^{s, \infty}(X)$ denotes the Besov space of order $s\in \mathbb{R}$.
\end{thm*}
We refer to this notion as the \underline{iterated regularity} characterization of conormal distributions.  A similar characterization exists for Lagrangian distributions \cite{Ho09}.  For more on Besov spaces, see \cite{Ste70}

Now, let ($T^*(X) \backslash 0,\omega_{T^*(X)}$), ($T^*(Y) \backslash 0,\omega_{T^*(Y)}$) be the cotangent bundles of the smooth manifolds $X$ and $Y$ with the 0-sections deleted and respective symplectic forms. ($T^*(X)\backslash 0 \times T^*(Y) \backslash 0, \omega_{T^*(X) \times T^*(Y)})$ is a symplectic manifold with respect to the twisted 2-form $\omega_{T^*(X) \times T^*(Y)} = \pi_L^*\omega_{T^*(X)} - \pi_R^*\omega_{T^*(Y)}$, where $\pi_L$ and $\pi_R$ are the left and right projections from our product space.  If $\Lambda \subset T^(X) \backslash 0 \times T^*(Y) \backslash 0$ is a conic Lagrangian submanifold with respect to this 2-form, then $\Lambda ' = \{\left( (x;\xi),(y ; \eta) \right) : \left( (x;\xi),(y ; -\eta) \right) \in \Lambda\}$ is called a \underline{canonical relation}; $\Lambda '$ is Lagrangian with respect to the twisted symplectic form $\pi_L^*\omega_{T^*(X)} - \pi_R^*\omega_{T^*(Y)}$. 

\begin{defn}
$F \in I^{m}(X,Y; C)$ is a \underline{Fourier integral operator of order m} (abbreviated by FIO) if the Schwartz kernel of $F$, $K_F(x,y)$, is an element of the space $I^{m}(X,Y; C')$.
\end{defn}

A natural question to consider: when is the composition of two FIOs again an FIO?  

\begin{thm} (Hormander 1971) \label{transintcalc}
Let $F_1 \in I^{m_1}(X,Y; C_1)$ and $F_2 \in I^{m_2}(Y,Z; C_2)$.  Suppose that $C_1' \times C_2' \pitchfork T^*(X) \times \triangle_{T^*(Y)} \times T^*(Z)$, where $\pitchfork$ denotes transverse intersection.  Then $F_1 \circ F_2 \in I^{m_1 + m_2}(X,Z;C_1 \circ C_2)$.
\end{thm}
\noindent Operator compositions that fall under the hypothesis of this theorem are said to satisfy the transverse intersection calculus.

\subsection{Paired lagrangians and nested conormals}

Classes of distributions associated with two cleanly intersecting Lagrangians manifolds were introduced by Melrose and Uhlmann \cite{MelUhl79} and Guillemin and Uhlmann \cite{GuiUhl81} in order to construct parametrices for systems of pseudodifferential operators that arose in various settings.  For the purposes of this paper, we define the class through multiphase functions \cite{Mend82} and symbol-valued symbols \cite{GU90}, as follows.

\begin{defn} \cite{Bott56}
Let $M$ and $N$ be smooth submanifolds of $X$.  Suppose that $M \cap N$ is also smooth.  Then for $p \in M \cap N$, $M$ and $N$ are said to be \underline{cleanly intersecting at $p$} if 
\begin{equation*}
T_p(M) \cap T_p(N) = T_p(M \cap N).
\end{equation*}
Moreover,  $(M,N)$ is a \underline{cleanly intersecting pair in codimension $k$} if $M$ and $N$ are cleanly intersecting for all $p \in M \cap N$ and $M \cap N$ is codimension $k$ in both submanifolds.
\end{defn}

\begin{defn}
Let $(\Lambda_0, \Lambda_1)$ be a cleanly intersecting pair of Lagrangians in codimension $k$ in $T^*(X) \backslash 0$.  Suppose $\lambda_0 \in \Lambda_0 \cap \Lambda_1$ and $\Gamma \subset X \times (\mathbb{R}^N \backslash 0) \times \mathbb{R}^k$ is an open conic set.  A \underline{multiphase function $\phi$} parametrizing the pair $(\Lambda_0, \Lambda_1)$ is a function $\phi(x; \theta; \sigma) \in C^{\infty}(\Gamma)$ such that 
\begin{enumerate}
\item $\phi_0(x; \theta) := \phi(x; \theta; 0)$ is a nondegenerate phase function parametrizing $\Lambda_0$ in a conic neighborhood of $\lambda_0$, and 
\item $\phi_1(x; (\theta, \sigma)) := \phi(x; \theta; \sigma)$ is a nondegenerate phase function parametrizing $\Lambda_1$ in a conic neighborhood of $\lambda_0$.
\end{enumerate}
\end{defn}

\noindent \underline{Example}:  The work of Guillemin and Uhlmann \cite{GuiUhl81} proves that any two pairs of cleanly intersecting Lagrangians are microlocally equivalent. One can thus consider the model pair $(\Lambda_0, \Lambda_1)$ in $T^*(\mathbb{R}^n)$ where $\Lambda_0$ and $\Lambda_1$ are the conormal bundles of $\{x=(x_1,...,x_n)=0\}$ and $\{x'=(x_{k+1},...,x_n)=0\}$, respectively, so that 
\begin{eqnarray*}
\Lambda_0 &=& T^*_0(\mathbb{R}^n) = \{(0; \xi) : \xi \in \mathbb{R}^n \backslash 0\}\\
\Lambda_1 &=& \{(x, \xi) : x_{k+1}=...=x_n=0, \xi_1=...=\xi_k=0\}.
\end{eqnarray*}
Then $\varphi(x; \theta ', \sigma) = x \cdot (\theta ', \sigma)$, for $(\theta ', \sigma) \in (\mathbb{R}^{n-k} \backslash 0) \times \mathbb{R}^k$ is an example of a multiphase parametrizing $(\Lambda_0, \Lambda_1)$.

The singularities of distributions represented by an oscillatory integral is dependent on the interaction between its phase and symbol.  As this discussion leads to a more general interpretation of Lagrangians distributions, it is natural to expect that the symbols themselves will generalize.
\begin{defn}
The space $S^{M_1, M_2}(X \times (\mathbb{R}^k_1 \backslash 0) \times \mathbb{R}^{k_2})$ of \underline{symbol-valued symbols} is the set of functions $a(x; \theta; \sigma) \in C^{\infty}(X \times \mathbb{R}^{k_1}  \times \mathbb{R}^{k_2})$ such that, for every $K \Subset X$, $(\alpha, \beta, \gamma) \in \mathbb{Z}_+^{k_1} \times \mathbb{Z}_+^{k_2} \times \mathbb{Z}_+^n$, the following differential estimate holds:
\begin{equation*}
|\partial_x^{\gamma} \partial_{\sigma}^{\beta }\partial_{\theta}^{\alpha} a(x; \theta; \sigma)| \lesssim  \langle \theta, \sigma \rangle ^{p - |\alpha|} \langle \sigma \rangle ^{l - |\beta|},
\end{equation*}
for $x \in K$.  If $|\theta| \gtrsim |\sigma|$ on the support of $a$, we say that $\theta$ is \underline{elliptic} to $\sigma$.
\end{defn}

Next, we define a generalized class of Fourier integral distributions associated with a cleanly intersecting pair of Lagrangians.  The following definition is a modification of the original formulation found in \cite{Mend82}.

\begin{defn}
Let $(\Lambda_0, \Lambda_1)$ be a cleanly intersecting pair of Lagrangians in codimension $k$ in $T^*(X) \backslash 0$.  The space of \underline{paired Lagrangian} distributions of order $p,l \in \mathbb{R}$ associated to $(\Lambda_0, \Lambda_1)$, denoted by $I^{p,l}(\Lambda_0, \Lambda_1)$, is the set of all locally finite sums of elements of $I^{p+1}(\Lambda_0) + I^{p}(\Lambda_1)$ and distributions of the form
\begin{equation}
u(x) = \int e^{i \phi(x; \theta; \sigma)} a(x; \theta; \sigma) d \theta d \sigma,
\end{equation}
where $a \in S^{\tilde{p}, \tilde{l}}(X \times (\mathbb{R}^N \backslash 0) \times \mathbb{R}^k)$, with $p = \tilde{p} + \tilde{l} + \frac{N+k}{2} - \frac{\dim X}{4}$, $l = - \tilde{l} - \frac{k}{2}$, and $\phi(x; \theta; \sigma)$ is multiphase parametrizing $(\Lambda_0, \Lambda_1)$ on a conic neighborhood of a point $\lambda_0 \in \Lambda_0 \cap \Lambda_1$.
\end{defn}

\noindent \underline{Example}: Consider $u(x) = H(x_1) \cdot \delta(x'')$ where $(x_1,x'') \in \mathbb{R}^n$, $H(x_1)$ is the Heaviside function in $x_1$, and $\delta(x'')$ is the delta function.  It follows that $u \in I^{\frac{n}{4}-\frac{1}{2}, -\frac{n}{4}-\frac{1}{2}}(N^*(\{x_1 \geq 0 \}), N^*(\{x''=0\}))$.  

The potential in our scattering problem will be a similar kind of paired Lagrangian distribution.

\begin{defn}
If $S_2 \subset S_1 \subset X$ are smooth submanifolds with codim$(S_1)=d_1$ and codim$(S_2)=d_1+d_2$, then $N^*(S_1)$ and $N^*(S_2)$ intersect cleanly in codimension $d_2$.  The space of distributions on $X$ conormal to the pair $(S_1, S_2)$, referred to as \underline{nested conormal distributions}, is
\begin{equation}
I^{M_1, M_2}(S_1, S_2) = I^{M_1+M_2+\frac{d_1+d_2}{2}-\frac{n}{4},-\frac{d_2}{2} - M_2}(N^*(S_1),N^*(S_2)).
\end{equation}
If $S_1$ and $S_2$ have defining functions $\{h_i(x)\}_{i=1}^{d_1}$ and $\{h_i(x)\}_{i=1}^{d_1+d_2}$, respectively, then $u \in I^{M_1, M_2}(S_1, S_2)$ can be locally expressed as 
\begin{equation}
\int_{\mathbb{R}^{d_1+d_2}} e^{i[\Sigma_i h_i(x) \theta_i]}a(x; \theta ', \theta '') d \theta ' d \theta '',
\end{equation}
where $a \in S^{M_1, M_2} (X \times (\mathbb{R}^{d_1} \backslash 0) \times \mathbb{R}^{d_2}).$
\end{defn}

Similar to the iterated regularity characterization of conormal distributions, there exists one for the nested conormals.
\begin{defn}
Suppose that $S_1$ and $S_2$ are submanifolds of $X$ that intersect cleanly.  Then $\mathcal{V}(S_1,S_2)$ is the space of vector fields that consists of those that are tangent to both $S_1$ and $S_2$.
\end{defn}

 \begin{thm*} \cite{GU90}
The space of nested conormal distributions with respect to $S_1, S_2$, for $S_2 \subset S_1$, denoted by $I^{M_1, M_2}(S_1, S_2)$, coincides with the set of all distributions $u \in \mathcal{D}'(X)$ such that for some $s_0$ depending on $M_1$ and $M_2$, for all $k \geq 0$,
\begin{equation*}
V_1...V_ku \in H^{s_0}_{loc}(X) 
\end{equation*}
with $V_j \in \mathcal{V}(S_1, S_2)$ for $1 \leq j \leq k$.  Here, $H^s_{loc}(\mathbb{R}^n)$ denotes the localized Sobolev space \cite{Foll95}.
\end{thm*}

\section{Nested triple conormals} \label{ntc}

In this section, we introduce a new class of distributions analogous to the nested pair conormals but with three Lagrangians intersecting pairwise cleanly with a smooth triple intersection.  The work of Mendoza and Uhlmann \cite{MendUhl84} develops a symbol calculus for a class of distributions associated to an intersecting system of three Lagrangians.  However the conditions imposed on their intersecting system are not satisfied in our situation, requiring us to take another approach.

\subsection{Iterated regularity characterization}
\begin{defn}
Suppose that $S_1, S_2,$ and $S_3$ are submanifolds of $X$ that have smooth pairwise intersections, a smooth triple intersection, and are in general position; see \cite[p.83]{GolGui73}.  Then $\mathcal{V}(S_1,S_2,S_3)$ is the space of vector fields consisting of those tangent to $S_1$, $S_2$, and $S_3$.
\end{defn}
\begin{defn}
The space of triple product-type conormal distributions with respect to $S_1, S_2,$ and $S_3$, $I(X;S_1, S_2, S_3)$, is the set of all distributions $u \in \mathcal{D}'(X)$ such that for some $s_0$ and all $k \geq 0$,
\begin{equation} \label{triplereg}
V_1...V_ku \in H^{s_0}_{loc}(X) 
\end{equation}
with $V_j \in \mathcal{V}(S_1, S_2, S_3)$ for $1 \leq j \leq k$.
\end{defn}

We will assume $S_3 \subset S_2 \subset S_1$ with codim$(S_1)=d_1$, codim$(S_2)=d_1+d_2$, and codim$(S_3)=d_1+d_2+d_3$ .  It is possible to introduce local coordinates in a neighborhood of $S_3$ with
\begin{eqnarray} \label{localcoords}
\nonumber S_1 &=& \{x_1=...=x_{d_1}\}\\
\nonumber S_2 &=& \{x_1=...=x_{d_1+d_2}\}\\
S_3 &=& \{x_1=...=x_{d_1+d_2+d_3}\};
\end{eqnarray}
denote the points in $\mathbb{R}^n$ as $(x',x'',x''',x'''')$ with $x'=(x_1,...,x_{d_1})$, $x''=(x_{d_1+1},...,x_{d_1+d_2})$, $x'''=(x_{d_1+d_2+1},...,x_{d_1+d_2+d_3})$, and $x''''=(x_{d_1+d_2+d_3+1},...,x_{n})$.  Let $(\xi ', \xi '', \xi ''', \xi '''')$ be the dual variables in phase space.  The following is based on the iterated regularity characterization for nested pairs in \cite{Melnotes89}; see also \cite{GU90}.

\begin{prop} \label{iterreg3}
If $u \in \mathcal{D}'(\mathbb{R}^n)$ and $S_3 \subset S_2 \subset S_1 \subset \mathbb{R}^n$ are as in (\ref{localcoords}), then $u \in I(\mathbb{R}^n; S_1, S_2, S_3)$ if and only if there exists an $s_0 \in \mathbb{R}$ such that
\begin{equation} \label{simpledo}
D_{x'}^{\alpha}D_{x''}^{\beta}D_{x'''}^{\nu}D_{x''''}^{\gamma}((x')^{\rho}(x'')^{\delta}(x''')^{\sigma}u) \in H^{s_0}(\mathbb{R}^n),
\end{equation} 
for all multi-indices $\alpha, \beta, \eta, \gamma, \rho, \delta, \sigma$ where $|\rho| \geq |\alpha|, |\delta| + |\rho| \geq |\alpha| + |\beta|, |\delta| + |\rho| + |\sigma| \geq |\alpha| + |\beta| + |\eta|$.
\end{prop}

\noindent \underline{Proof}:  We find a local set of generators, over $C^{\infty}(\mathbb{R}^n)$, for the algebra $\mathcal{V}(S_1, S_2, S_3)$ consisting of the differential operators that appear in (\ref{triplereg}); we will do induction on the order of the operator in (\ref{simpledo}).

Take $V= \Sigma_{i=1}^{d_1} a_i(x)D_{x_i} + \Sigma_{j=d_1+1}^{d_1+d_2} b_j(x)D_{x_j} + \Sigma_{k=d_1+d_2+1}^{d_1+d_2+d_3} c_k(x)D_{x_k}$ \newline $+ \Sigma_{\ell=d_1+d_2+d_3+1}^{n}d_{\ell}(x)D_{x_{\ell}}$, with $a_i,b_j,c_k,d_{\ell} \in C^{\infty}(\mathbb{R}^n)$, then $V$ being tangent to $S_1$, $S_2$, and $S_3$ means that the $a_i=0$ at $x'=0$, $a_i=b_j=0$ at $x'=x''=0$, and $a_i=b_j=c_k=0$ at $x'=x''=x'''=0$, respectively.  Then, 
\begin{eqnarray} \label{localbasis}
\nonumber x_iD_{x_j} & \text{ for }& i \leq d_1, j \leq d_1,\\
\nonumber x_iD_{x_j} & \text{ for }& i \leq d_1+d_2, d_1 < j \leq d_1 + d_2,\\
\nonumber x_iD_{x_j} & \text{ for }& i \leq d_1+d_2+d_3, d_1 + d_2 < j \leq d_1 + d_2 + d_3, \text{ and }\\
D_{x_k} & \text{ for }& k > d_1 + d_2 + d_3,
\end{eqnarray}
gives our desired basis.

Assume (\ref{simpledo}) holds for $|\alpha| + |\beta| + |\eta| + |\gamma| \leq p$, where the induction step will be done on $p$.  The application of an operator $V_{p+1} \in \mathcal{V}(S_1, S_2, S_3)$, which is a linear combination of the vector fields in (\ref{localbasis}), to the expression in (\ref{simpledo}) gives 
\begin{eqnarray}
V_1...V_{p+1} &=& \Sigma a^{\rho \delta \sigma i j}_{\alpha \beta \eta \gamma} D_{x'}^{\alpha}D_{x''}^{\beta}D_{x'''}^{\nu}D_{x''''}^{\gamma}((x')^{\rho}(x'')^{\delta}(x''')^{\sigma}\phi_{i j}x_i D_{x_j})\\
&+& \Sigma b^{\rho \delta \sigma k}_{\alpha \beta \eta \gamma} D_{x'}^{\alpha}D_{x''}^{\beta}D_{x'''}^{\nu}D_{x''''}^{\gamma}((x')^{\rho}(x'')^{\delta}(x''')^{\sigma}\phi_{i j}\phi_k D_{x_k}),
\end{eqnarray}
where $i,j,k$ have the restrictions indicated in (\ref{localbasis}) and the coefficients are in $C^{\infty}$.  Further commutation of the differential operators with the coefficients $x', x'',x''', \phi_{ij},$ and $\phi_j$ gives the expression in (\ref{simpledo}) with $|\alpha| + |\beta| + |\eta| + |\gamma| \leq p+1$.  $\blacksquare$

\begin{prop}
If $u \in I(X;S_1, S_2, S_3)$, then $WF(u) \subseteq N^*(S_1) \cup N^*(S_2) \cup N^*(S_3)$.
\end{prop}

\noindent \underline{Proof}: Using Definition \ref{triplereg} and the vector fields computed in the proof of Proposition \ref{iterreg3}, it follows that $u$ is not singular in the $x''''$ variable.  $\blacksquare$

The iterated regularity definition is not limited to the configuration of nested submanifolds.  The computation of a basis for $\mathcal{V}(S_1, S_2, S_3)$ for a more general triple needs additional assumptions and is more complicated.  This is left to the interested reader.

\subsection{Oscillatory representation}

\begin{defn}
The space $S^{M_1, M_2, M_3}(X \times (\mathbb{R}^k_1 \backslash 0) \times \mathbb{R}^{k_2}  \times \mathbb{R}^{k_3})$ of \underline{triple symbol-} \underline{valued symbols} is the set of functions $a(x; \theta; \sigma; \tau) \in C^{\infty}(X \times \mathbb{R}^{k_1} \times \mathbb{R}^{k_2} \times \mathbb{R}^{k_3} )$ such that, for every $K \Subset X$, $(\alpha, \beta, \eta, \gamma) \in \mathbb{Z}_+^{k_1} \times \mathbb{Z}_+^{k_2} \times \mathbb{Z}_+^{k_3} \times \mathbb{Z}_+^n$, the following differential estimate holds:
\begin{equation*}
|\partial_x^{\gamma} \partial_{\tau}^{\eta }\partial_{\sigma}^{\beta } \partial_{\theta}^{\alpha} a(x; \theta; \sigma; \tau)| \lesssim  \langle \theta, \sigma, \tau \rangle ^{M_1 - |\alpha|} \langle \sigma, \tau \rangle ^{M_2 - |\beta|} \langle \tau \rangle ^{M_3 - |\eta|},
\end{equation*}
for $x\in K$.  
\end{defn}

\begin{prop}
Let $u \in I(\mathbb{R}^n; S_1, S_2, S_3)$ with the local coordinates in (\ref{localcoords}).  Then 
\begin{equation} \label{triposcill}
u = \int_{\mathbb{R}^{d_1+d_2+d_3}} e^{i(x'\cdot \xi ' + x''\cdot \xi '' + x'''\cdot \xi ''')} a(x; \xi '; \xi ''; \xi ''') d \xi ' d\xi '' d\xi ''',
\end{equation}
where $a \in S^{M_1, M_2, M_3}(\mathbb{R}^n \times (\mathbb{R}^d_1 \backslash 0) \times \mathbb{R}^{d_2}  \times \mathbb{R}^{d_3})$.
\end{prop}
\noindent \underline{Proof}:  Following Melrose's original line of argument, also found in \cite{GU90}, we can assume that $u$ is compactly supported; otherwise, we can multiply by a smooth cutoff supported near $S_3$.  Take the Fourier transform in $(x',x'',x''')$ and using (\ref{simpledo}) we get
\begin{eqnarray}
\nonumber &&(\xi ')^{\alpha}(\xi '')^{\beta}(\xi ''')^{\eta}D_{x''''}^{\gamma}D_{\xi '}^{\rho}D_{\xi ''}^{\delta}D_{\xi '''}^{\sigma}\\
\nonumber && \in L^2(\mathbb{R}^{d_1} \backslash 0 \times \mathbb{R}^{d_2} \backslash 0 \times \mathbb{R}^{d_3} \backslash 0 \times \mathbb{R}^{n-d_1-d_2-d_3}; \\
&& \langle \xi ', \xi '', \xi ''' \rangle^{s_0} d \xi ' d \xi '' d \xi ''' d x'''')
\end{eqnarray}
for some $s_0$ when $|\rho| \geq |\alpha|, |\delta| + |\rho| \geq |\alpha| + |\beta|, |\delta| + |\rho| + |\sigma| \geq |\alpha| + |\beta| + |\eta|$.  The Sobolev embedding theorem \cite{Foll95} tells us $a(x; \xi '; \xi ''; \xi ''') := \hat{u}(\xi ', \xi '', \xi ''', x '''')$ satisfies a triple symbol-valued symbol estimate for some $M_1, M_2,$ and $M_3$ depending on the dimension of our submanifolds $S_1$, $S_2$, and $S_3$.  $\blacksquare$

We will now generalize the proof of Proposition 1.20 in \cite{GU90} to show that distributions of the form (\ref{triposcill}) are independent of the choice of coordinates.  

Suppose $u$ is of the form (\ref{triposcill}) with an integral that is absolutely convergent; otherwise, integrate by parts to lower the order of the triple symbol-valued symbol.  A change of variables that preserves (\ref{localcoords}) must be of the form 
\begin{equation}
\begin{cases}
x_i = \sum_{j=1}^{d_1} A_{ij}(y), & 1 \leq i \leq d_1 \\
x_i = \sum_{j=1}^{d_1+d_2} B_{ij}(y), & d_1 < i \leq d_1 + d_2 \\
x_i = \sum_{j=1}^{d_1+d_2+d_3} C_{ij}(y), & d_1 +d_2 < i \leq d_1 + d_2 + d_3 \\

\end{cases}
\end{equation}
with $A_{ij}, B_{ij},$ and $C_{ij}$ being smooth.  Plugging this new coordinate transformation into (\ref{triposcill}) and setting 
\begin{equation}
\begin{cases}
\Xi_i = \sum_{j=1}^{d_1} A_{ji}\xi_j + \sum_{j=d_1+1}^{d_1+d_2} B_{ji}\xi_j + \sum_{j=d_1 + d_2+ 1}^{d_1+d_2+d_3} C_{ji}\xi_j, & 1 \leq i \leq d_1 \\
\Xi_i = \sum_{j=d_1+1}^{d_1+d_2} B_{ji}\xi_j + \sum_{j=d_1 + d_2+ 1}^{d_1+d_2+d_3} C_{ji}\xi_j, & d_1 < i \leq d_1 + d_2 \\
\Xi_i = \sum_{j=d_2+d_3+1}^{d_1+d_2+d_3} C_{ji}\xi_j, & d_1 +d_2 < i \leq d_1  \\
& + d_2 + d_3,
\end{cases}
\end{equation}
gives the new oscillatory representation of 
\begin{equation}
u(y) = \int e^{i(\Xi ' \cdot y' + \Xi '' \cdot y'' + \Xi ''' \cdot y''')} b(y; \Xi ' ; \Xi '' ; \Xi ''') d \Xi ' d \Xi '' d \Xi '''.
\end{equation}
Here, the interested reader can prove $b = a(x(y),\xi'(y,\Xi',\Xi'',\Xi'''),\xi''(y,\Xi',\Xi''$ $,\Xi'''), \xi  '''(y, \Xi ', \Xi '', \Xi ''')) \times \left | \frac{D \xi}{D \Xi} \right |$ is another triple symbol-valued symbol.  This shows that (\ref{triposcill}) is independent of the choice of coordinates.  We can now make the following definition.
\begin{defn} \label{tripinva}
Let $S_3 \subset S_2 \subset S_1$ be submanifolds of $X$ with codimensions $d_1 + d_2 + d_3, d_1 + d_2,$ and $d_1$, respectively.  Then $I^{M_1, M_2, M_3}(\mathbb{R}^n; S_1, S_2, S_3)$ is the space of locally finite sums of distributions of the form 
\begin{equation}
\int_{\mathbb{R}^{d_1+d_2+d_3}} e^{i(H_1(x) \cdot \xi ' + H_2(x) \cdot \xi '' + H_3(x)\cdot \xi ''')} a(x; \xi '; \xi ''; \xi ''') d \xi ' d\xi '' d\xi ''',
\end{equation}
where $a \in S^{M_1, M_2, M_3}(\mathbb{R}^n \times (\mathbb{R}^d_1 \backslash 0) \times \mathbb{R}^{d_2}  \times \mathbb{R}^{d_3})$ and $\{H_1(x)\}, \{H_1(x)=H_2(x)\},$ and $\{H_1(x)=H_2(x)=H_3(x)\}$ are defining functions for $S_1, S_2$, and $S_3$, respectively.
\end{defn}

\begin{prop} \label{microdub}
Let $u \in I^{M_1, M_2, M_3}(\mathbb{R}^n; S_1, S_2, S_3)$.  Then, away from $N^*(S_1) \cap N^*(S_2) \cap N^*(S_3)$, 
\begin{eqnarray}
\nonumber u &\in& I^{M_1+M_2+ \frac{d_1+d_2}{2} - \frac{n}{4}, -M_2 - \frac{d_2}{2}}(N^*(S_1),N^*(S_2)) \\
\nonumber &+& I^{M_1+M_2+M_3+ \frac{d_1+d_2+d_3}{2} - \frac{n}{4}, -M_2 -M_3 - \frac{d_2+d_3}{2}}(N^*(S_1),N^*(S_3))\\ 
&+& I^{M_1+M_2+M_3 + \frac{d_1+d_2+d_3}{2} - \frac{n}{4}, -M_3 - \frac{d_3}{2}}(N^*(S_2),N^*(S_3))
\end{eqnarray}
\end{prop}

\noindent \underline{Proof}: By using Definition \ref{tripinva}, we can assume 
\begin{equation}
u = \int_{\mathbb{R}^{d_1+d_2+d_3}} e^{i(x'\cdot \xi ' + x''\cdot \xi '' + x'''\cdot \xi ''')} a(x; \xi ', \xi '', \xi ''') d \xi ' d\xi '' d\xi '''
\end{equation}
for $a \in S^{M_1,M_2,M_3}$.  Therefore, 
\begin{eqnarray}
\nonumber N^*(S_1) &=& \{(0,x'',x''',x'''';\xi ',0,0,0) : \xi ' \neq 0\}, \\
\nonumber N^*(S_2) &=& \{(0,0,x''',x'''';\xi ',\xi '',0,0) : (\xi ', \xi '') \neq 0\}, \text{ and }\\
N^*(S_3) &=& \{(0,0,0,x'''';\xi ',\xi '',\xi ''',0) : (\xi ', \xi '', \xi ''') \neq 0\}
\end{eqnarray}
which shows
\begin{equation} \label{tripconint}
N^*(S_1) \cap N^*(S_2) \cap N^*(S_3) = \{(0,0,0,x'''';\xi ',0,0,0): \xi ' \neq 0\}.
\end{equation}
Away from (\ref{tripconint}), we have $(x',x'',x''') \neq 0$, $\xi '' \neq 0$, or $\xi ''' \neq 0$.  In the first case, integration by parts in the $\xi', \xi '',$ or $\xi '''$ variables contributes a $C^{\infty}$ function, an element of $I^{M_1}(S_1)$, or an element in $I^{M_1, M_2}(S_1,S_2)$. 

Now, suppose $\xi '' \neq 0$.  Then either $\{|\xi ''| \geq |\xi '|\}$ or $\{|\xi ''| \geq |\xi '''|\}$.  Assuming the former case, we get that $u \in I^{M_1+M_2+M_3 + \frac{d_1+d_2+d_3}{2} - \frac{n}{4}, -M_3 - \frac{d_3}{2}}(N^*(S_2),$ $N^*(S_3))$ since $a \in S^{M_1+M_2, M_3}$.  Due to being away from the intersection described in (\ref{tripconint}), $\xi '$ is allowed to be 0.  As $\xi '$ cannot be elliptic to all the other variables, for otherwise we would be localized to the intersection in (\ref{tripconint}), it follows that we are either again in the region $\{|\xi ''| \geq |\xi '|\}$ or in the region $\{|\xi '''| \geq |\xi '|\}$; the second region is discussed next.  

Consider $\xi ''' \neq 0$.  Then either $\{|\xi '''| \geq |\xi '|\}$ or $\{|\xi '''| \geq |\xi ''|\}$.  The former case again shows that 
\begin{equation}
u \in I^{M_1+M_2+M_3 + \frac{d_1+d_2+d_3}{2} - \frac{n}{4}, -M_3 - \frac{d_3}{2}}(N^*(S_2),N^*(S_3)),
\end{equation}
while the latter case implies  
\begin{equation}
u \in I^{M_1+M_2+M_3+ \frac{d_1+d_2+d_3}{2} - \frac{n}{4}, -M_2 -M_3 - \frac{d_2+d_3}{2}}(N^*(S_1),N^*(S_3)).   
\end{equation}
as $a \in S^{M_1, M_2 + M_3}$.  Each of the paired Lagrangian distributions found falls into a class described in the statement of the proposition.  $\blacksquare$

\section{Lagrangian submanifolds and canonical relations} \label{lscr}

\noindent Each of the subsections below will be focused on the specific geometry induced by a given operator.  In each case, we give the corresponding canonical relations for the operator, analyze the compositions with certain Lagrangian manifolds, and state properties of the resulting manifolds.

\subsection{Geometry determined by a multiplication operator}

\noindent Recall for $q \in C^{\infty}_0(\mathbb{R}^n)$, multiplication by $q$, which we denote by $M_q$, is a pseudodifferential operator of order 0.  Moreover, $M_q$ for $q$ conormal of order $\mu$ is an example of, in the language of \cite{GU90}, a pseudodifferential operator with a singular symbol of order $0$ and $\mu$ when using conormal notation.  This notation will be described further in a later section.

For the operator $M_q$ with $q \in I^{M_1,M_2}(S_1, S_2)$, there are three canonical relations in $T^*(\mathbb{R}^n \times \mathbb{R} \times S^{n-1}) \times T^*(\mathbb{R}^n \times \mathbb{R} \times S^{n-1})$ associated to it:
\begin{eqnarray}
\triangle_{T^*(\mathbb{R}^n \times \mathbb{R} \times S^{n-1}) \times T^*(\mathbb{R}^n \times \mathbb{R} \times S^{n-1})},
\end{eqnarray}
\begin{eqnarray} \label{can rel 1}
 C_{S_1} & = & \{ \left( (x,t,\omega ; \xi, \tau, \Omega), (x,t,\omega ; \xi + D_x((h_i)_{i=1}^{d_1})\cdot \theta, \tau, \Omega) \right): \\
\nonumber & & (x,t,\omega ; \xi, \tau, \Omega) \in T^*(\mathbb{R}^n \times \mathbb{R} \times S^{n-1})_{|S_i}, \text{ and } \theta \in \mathbb{R}^{d_1} \backslash 0\},
\end{eqnarray}
where $\{h_i\}_{i=1}^{d_1} \text{ defines } S_1$ and 
\begin{eqnarray} \label{can rel 2}
C_{S_2} & = & \{ \left( (x,t,\omega ; \xi, \tau, \Omega), (x,t,\omega ; \xi + D_x((h_i)_{i=1}^{d_1+d_2})\cdot \theta, \tau, \Omega) \right): \\
\nonumber & & (x,t,\omega ; \xi, \tau, \Omega) \in T^*(\mathbb{R}^n \times \mathbb{R} \times S^{n-1})_{|S_2}, \text{ and } \theta \in \mathbb{R}^{d_1+d_2} \backslash 0\},
\end{eqnarray}
where $\{h_i\}_{i=1}^{d_1+d_2} \text{ defines } S_2$.
For our purposes, the three most important Lagrangians in $T^*(\mathbb{R}^n \times \mathbb{R} \times S^{n-1})$ that interact with the canonical relations of $M_q$ are
\begin{eqnarray}
\Lambda_+ &=& \{ (x,x \cdot \omega, \omega ; -\sigma \omega, \sigma, -\sigma i^*_{\omega}(x)) : \\
\nonumber && (x,t, \omega) \in \mathbb{R}^n \times \mathbb{R} \times S^{n-1}, \sigma \in \mathbb{R} \backslash 0 \} \\
\Lambda_1 &=& \{ (y, t, \omega ; \nu , 0, 0) : (y,\nu) \in N^*(S_1) \} \\
\Lambda_2 &=& \{ (y, t, \omega ; \nu , 0, 0) : (y,\nu) \in N^*(S_2) \}.
\end{eqnarray}
It is important to point of out these three Lagrangians are also the conormal bundles of $S_+, S_1, \text{ and } S_2$, respectively, where $S_+=\{t - x \cdot \omega = 0: (x,t,\omega) \in \mathbb{R}^n \times \mathbb{R} \times S^{n-1} \}$.  As $S_i \pitchfork S_+$ for $i=1,2$, it follows $S_i \cap S_+ = S_{i+}$ are smooth submanifolds of $\mathbb{R}^n \times \mathbb{R} \times S^{n-1}$.  Composition of the canonical relations in (\ref{can rel 1}) and (\ref{can rel 2}) with $\Lambda_+$ give  
\begin{eqnarray} 
\label{lambda1+} \Lambda_{1+} &=& \{ (y,y \cdot \omega, \omega ; \nu -\sigma \omega, \sigma, -\sigma i^*_{\omega}(x)) : \\
\nonumber && (y, \nu) \in N^*(S_1), \omega \in S^{n-1}, \sigma \in \mathbb{R} \backslash 0 \},\text{ and } \\
\label{lambda2+} \Lambda_{2+} &=& \{ (y,y \cdot \omega, \omega ; \nu -\sigma \omega, \sigma, -\sigma i^*_{\omega}(x)) : \\
\nonumber && (y, \nu) \in N^*(S_2), \omega \in S^{n-1}, \sigma \in \mathbb{R} \backslash 0 \}.
\end{eqnarray}
Once again the Lagrangians submanifolds described in (\ref{lambda1+}) and (\ref{lambda2+}) are the conormal bundles of $S_{1+}$ and $S_{2+}$.  From this fact, it follows immediately that $(\Lambda_{1+},\Lambda_{2+})$ is a cleanly intersecting pair in codimension $d_2$.  Similar reasoning can be applied to other pairs of the above Lagrangians.

\subsection{Geometry determined by the wave equation}

\noindent The operator $\square^{-1}$, to be elaborated upon in Section 7, is another example pseudodifferential operator with a singular symbol.  More, specifically $\square^{-1} \in I^{-\frac{3}{2}, -\frac{1}{2}}(\triangle_{T^*(\mathbb{R}^n \times \mathbb{R} \times S^{n-1}) \times T^*(\mathbb{R}^n \times \mathbb{R} \times S^{n-1})}, C_{\square})$, where 
\begin{eqnarray}
\nonumber C_{\square} &=& \{ \left((x,t, \omega; \xi, |\xi|, \Omega), (x + (t-s)\frac{\xi}{|\xi|}, s, \omega; \xi, | \xi |, \Omega) \right): \\
\nonumber & & (x,t,\omega) \in \mathbb{R}^n \times \mathbb{R} \times S^{n-1}, s \in \mathbb{R} - 0, \xi \in \mathbb{R}^n, \\ 
&& \Omega \in T_{\omega}^*(S^{n-1}) \}
\end{eqnarray}
is the flowout relation generated by the characteristic variety 
\begin{equation}
char(\square) \break = \{(x,t, \omega ; \xi, \tau, \Omega) : |\tau|^2 = |\xi|^2 \}.
\end{equation} 

The application of this canonical relation to certain Lagrangians yields new Lagrangians away from a prescribed set.  As $\Lambda_i$, for $i=1,2$, does not meet $char(\square)$, it follows that $C_{\square} \circ \Lambda_i = \emptyset$ and any application of $\square^{-1}$ to distributions with $WF$ not supported there will give a $C^{\infty}$ function.  Also, since $\Lambda_+ \subset char(\square)$, $C_{\square} \circ \Lambda_+ = \Lambda_+$.

Let 
\begin{eqnarray} \label{char1}
\nonumber \Sigma^1 &=& \Lambda_{1+} \cap char(\square) \\
&=& \{|\nu_1|^2 - 2\nu_1 \cdot \omega = 0: (y_1,\nu_1) \in N^*(S_1), \omega \in S^{n-1} \}
\end{eqnarray}
and 
\begin{eqnarray} \label{char2}
\nonumber \Sigma^2 &=& \Lambda_{2+} \cap char(\square) \\
 &=&\{|\nu_2|^2 - 2\nu_2 \cdot \omega = 0: (y_2,\nu_2) \in N^*(S_2), \omega \in S^{n-1}\}.
\end{eqnarray}
It is from this set that our new smooth Lagrangians will be flowed out.  However, there exist points for which these flowouts are not guaranteed to be smooth, therefore demanding that we impose some further restrictions.

From (\ref{char1}) and (\ref{char2}), if $\nu_1 \cdot \omega = 0$ or $\nu_2 \cdot \omega = 0$, then $\nu_1$ or $\nu_2$ is forced to be 0 therefore putting $0$-sections into $\Lambda_1$ or $\Lambda_2$, which violates the 0-section assumption.  Notice that this dot product is 0 when $\omega$ is tangent to the submanifolds $S_{1+}$ or $S_{2+}$ (In \cite{GU93}, such elements of $S^{n-1}$ are referred to as \underline{tangential rays}).  Hence, it is natural to consider the open region where  $\omega$ is not tangent to $S_1$ or $S_2$.  We define 
\begin{eqnarray*} 
\Sigma_1 := \Sigma^1 \cap \{\nu_1 \cdot \omega \neq 0: (y_1,\nu_1) \in N^*(S_1), \omega \in S^{n-1} \}
\end{eqnarray*} 
and 
\begin{eqnarray*}
\Sigma_2 := \Sigma^2 \cap \{\nu_2 \cdot \omega \neq 0: (y_2,\nu_2) \in N^*(S_2), \omega \in S^{n-1} \}.
\end{eqnarray*} 
After solving our characteristic equation, we arrive at 
\begin{eqnarray}
\label{good set 1} \Sigma_1 &=& \{ (y_1,y_1 \cdot \omega, \omega ; \nu_1 -\sigma \omega, \sigma, -\sigma i^*_{\omega}(x)) : (y_1, \nu_1) \in N^*(S_1), \\
\nonumber && \omega \in S^{n-1}, \sigma = \frac{|\nu_1|^2}{2 \nu_1 \cdot \omega} \},  \text{ and }\\
 \label{good set 2} \Sigma_2 &=& \{ (y_i,y_i \cdot \omega, \omega ; \nu_2 -\sigma \omega, \sigma, -\sigma i^*_{\omega}(x)) : (y_2, \nu_2) \in N^*(S_2), \\
\nonumber &&  \omega \in S^{n-1}, \sigma = \frac{|\nu_2|^2}{2 \nu_2 \cdot \omega} \}.
\end{eqnarray}
This reparametrization is done to show the above sets are smooth submanifolds of $\Lambda_{1+}$ and $\Lambda_{2+}$, respectively.  

The Hamiltonian vector field associated to the symbol of $\square$ is 
\begin{eqnarray*}
H_{\square} = -\xi \cdot \frac{\partial}{\partial x} + \tau \frac{\partial}{\partial t},
\end{eqnarray*} 
where $x$ and $t$ are the spatial and time coordinates on $T^*(\mathbb{R}^n \times \mathbb{R} \times S^{n-1})$.  In order to guarantee our flowouts are smooth Lagrangians, we must start at  points of $\Sigma_1$ or $\Sigma_2$ at which $H_{\square}$ is transverse.  In order for this to happen, $H_{\square}$ cannot be in the tangent space of $\Sigma_i$, for $i=1,2$.  For if it were, then $-\xi \cdot \frac{\partial}{\partial x}$ will lie in the horizontal tangent space of the manifolds parametrized in (\ref{good set 1}) and (\ref{good set 2}).  But under these parametrizations, $-\xi = -(\nu - \sigma \omega)$ and $\tau = \sigma$.  The only way this can occur, in $\Sigma_1$ for example, is if $\nu = 0$ and $\omega$ lies in the tangent space of $S_{1+}$, with either violating the conditions in (\ref{good set 1}).  Hence $H_{\square} \pitchfork \Sigma_i$ for $i=1,2$ and generates $H_{\square}$ the following Lagrangian manifolds:
\begin{eqnarray}
\label{flowout1} \Lambda_{1+}^a &=& \{ (y - r(\nu -\sigma \omega),y \cdot \omega + r, \omega ; \nu -\sigma \omega, \sigma, -\sigma i^*_{\omega}(y)) : \\
\nonumber & & (y, \nu) \in N^*(S_1), \omega \in S^{n-1}, r \in \mathbb{R}-0, \sigma = \frac{|\nu|^2|}{2 \nu \cdot \omega} \}\\
\label{flowout2} \Lambda_{2+}^a &=& \{ (y - r(\nu -\sigma \omega),y \cdot \omega + r, \omega ; \nu -\sigma \omega, \sigma, -\sigma i^*_{\omega}(y)) : \\
\nonumber & & (y, \nu) \in N^*(S_2), \omega \in S^{n-1}, r \in \mathbb{R}-0, \sigma = \frac{|\nu|^2|}{2 \nu \cdot \omega} \}
\end{eqnarray}

\vskip 1cm

\begin{figure} \label{schemlag}
\begin{center}
\begin{pspicture}(-5,0)(5,6)
  \psline(0,0)(-2,5)
  \psline(0,0)(2,5)
  \psellipse(0,5)(2,0.5)
  \psline(0,0)(0,4.5)
  \psbezier(0,0)(0.1,1.1)(0,2)(1,4.59)   
  \psbezier(0,0)(-0.05,1.1)(-0,2)(-1,4.59)   
  \psbezier(-5,2)(-4,4)(-1.2,2)(-0.2,2)
  \psbezier(-4,1)(-3,3)(-2,1.8)(-0.8,1.8)
  \psbezier(-0.8,1.8)(-0.5,1.8)(-0.5,2)(-0.2,2)
  \psbezier(-0.2,2)(0.1,2)(0.1,2.2)(0.4,2.2)
  \psbezier(5.2,2.2)(4.2,4.2)(1.4,2.2)(0.4,2.2)
  \psbezier(4.2,1.2)(3.2,3.2)(2.2,2)(1,2)
  \psbezier(1,2)(0.7,2)(0.7,2.2)(0.4,2.2)
  \rput{0}(5.2,3){$\Lambda_{2+}$}
  \rput{0}(4.2,2){$\Lambda_{1+}$}
  \rput{0}(0.4,3.5){$\Lambda_{+}$}
  \rput{0}(-.8,3.5){$\Lambda_{1+}^a$}
  \rput{0}(1.2,4){$\Lambda_{2+}^c$}
\end{pspicture}
\end{center}
\begin{center}
Figure \ref{schemlag}: A schematic representation of the Lagrangians; the cone represents $char(\square)$.
\end{center}
\end{figure}  

As all of our parametrizations come from a single coordinate patch, we can set coordinates equal to one another, solve the resulting equations, and find parametrizations for the intersections.  Doing this for $\Lambda_+$ and $\Lambda_{1+}$ shows $\nu = 0$.  Similarly, when looking at the parametrization of $\Lambda^a_{1+}$, evaluating the expression $\frac{|\nu|^2}{2 \nu \cdot \omega}$ at $\nu = 0$ makes $\sigma=0$, violating the 0-section assumption once more.  This shows $\Lambda_+ \cap \Lambda_{1+} \subset \Sigma^1 \backslash \Sigma_1$.  The exact same reasoning gives $\Lambda_+ \cap \Lambda_{1+} \subset \Sigma^2 \backslash \Sigma_2$.

The differential of the projection $\pi$ from the $\Lambda^a_{1+}$ onto the spatial variables $(x,t, \omega)$ is $\frac{D(x,t, \omega)}{D(y,r,\theta, \omega)} =$
\begin{eqnarray}
\begin{pmatrix}
  j^*(I_{n-d_1} + rd_y(\nu -\sigma \omega)) & \nu - \sigma \omega & r \kappa^*d_{\theta}(\nu - \sigma \omega) & ri^*_{\omega}(d_{\omega}(\sigma\omega))\\
  j^*(\omega - rd_y(\sigma))  & -\sigma & -r \kappa^*d_{\theta}(\sigma) & i^*_{\omega}(y)\\
  0 & 0 & 0 & I_{n-1}
\end{pmatrix}
\end{eqnarray}
where $j: T_{y}(S_{1}) \hookrightarrow T_{y}(\mathbb{R}^n)$, $i: T_{\omega}(S^{n-1}) \hookrightarrow T_{\omega}(\mathbb{R}^n)$, and $\kappa: d_1 \text{-plane} \hookrightarrow T_{y}(\mathbb{R}^n)$ are the differentials of the inclusions of the respective submanifolds.  By using Euler's identity for homogeneous functions, it follows $\frac{\partial}{\partial r} \in \text{span}$ $\{\frac{\partial}{\partial \theta_j}\}_{j=1}^{d_1}$.  As $y$ is given by $n-d_1$ parameters and $\theta$ by $d_1$, we see that the rank of this projection is locally constant and equal to $2n-1$ if the rank of the submatrix 
\begin{equation} \label{hypersub}
\begin{pmatrix}
j^*(I_{n-d_1} + rd_y(\nu -\sigma \omega)) \\
j^*(\omega - rd_y(\sigma)) 
\end{pmatrix}  
\end{equation}
is locally constant and equal to $n-d$.  This will hold away from the lower-dimensional set $\{(-r)^{n-d_1}\det(-d_y(\nu -\sigma \omega)- \frac{1}{r}I_{n-d_1}) = 0\}$, which is the singular set for the top submatrix in (\ref{hypersub}).  We remind ourselves that $r \neq 0$ for otherwise our parametrization drops rank.  Notice this determinant gives a polynomial in $r$, having only finitely many solutions in this parameter.  The same arguments are used to show that the rank of $\pi$ from $\Lambda_{2+}^c$ is also $2n-1$.  We use \cite[Prop. 3.7.2]{Duis11} on conormal bundles to help conclude 

\begin{lem}
Away from a lower-dimensional set in phase space, $\Lambda_{1+}^a$ and $\Lambda_{2+}^c$ are the conormal bundles of hypersurfaces in $\mathbb{R}^n \times \mathbb{R} \times S^{n-1}$.
\end{lem} 
\noindent This fact will play a substantial role in the approach we take when solving our desired inverse problem, particularly in Section 8. 

We will now recall a simple fact from linear algebra.  Suppose that $W_1$ and $W_2$ are subspaces of a finite dimensional vector space $V$.  We know $\dim(W_1 + W_2) = \dim(W_1) + \dim(W_2) - \dim(W_1 \cap W_2)$.  If we consider two submanifolds $M$ and $N$ of $X$, and assume $p \in M \cap N$ with the intersection being smooth.  As $T_p(M \cap N) \subset T_p(M) \cap T_p(N)$, it follows $$\dim(T_pM + T_pN) \leq \dim(T_pM) + \dim(T_pN) - \dim(T_p(M \cap N)).$$  For the given flowouts, an upcoming calculation shows $\Lambda_{1+}^a \cap \Lambda_{2+}^c$ has dimension $2n-d_2$.  For $\lambda \in \Lambda_{1+}^a \cap \Lambda_{2+}^c$, $\dim(T_{\lambda}\Lambda_{1+}^a + T_{\lambda}\Lambda_{2+}^a)$ cannot be greater than $2n+d_2$. Hence, if we can prove this sum of tangent spaces has dimension at least $2n+d_2$, then $\dim(T_p(M \cap N)) = \dim(T_p(M) \cap T_p(N))$ and the intersection is clean. 

Let us label each parameter that appears in (\ref{flowout2}) with the symbol $\text{ } \bar{ } \text{ }$ above it. Also, let $\theta ' \in \mathbb{R}^{d_1} \backslash 0$ parametrize the fibers of $\Lambda_{1+}^a$ and $(\bar{\theta}', \bar{\theta}'') \in \mathbb{R}^{d_1+d_2} \backslash 0$ parametrize those in $\Lambda_{2+}^c$.  Observing the spatial spherical and $\tau$ coordinates imply that $\sigma = \bar{\sigma}$ and $\omega = \bar{\omega}$, the equality in the $\xi$-coordinates tells us that $\nu=\bar{\nu}$ and that the fibers which are characteristic in $N^*(S_{1+})$ must coincide with those that are characteristic in $N^*(S_{2+})$.  Furthermore, after setting the $x$ coordinates equal and using the facts above, $y = \bar{y}$.  It is worth noting $\sigma=\bar{\sigma}$ is actually an implication of $N^*(S_{1+})$ intersecting $N^*(S_{2+})$.  In conclusion, it follows that in either parametrization of our flowouts, the intersection is obtained by simply restricting ourselves to $S_{2+}$ in the $(x,t)$ coordinate or to the fibers of $N^*(S_{1+})$ in the $(\xi, \tau)$ coordinates in the case of $\Lambda_{1+}^a$ or $\Lambda_{2+}^c$, respectively.  As this restriction only drops rank by $d_2$, $\Lambda_{1+}^a \cap \Lambda_{2+}^c$ has codimension $d_2$ in both flowouts.

To show the dimension of the tangents do add up correctly, we concatenate the differentials of the parametrizations in (\ref{flowout1}) and (\ref{flowout2}) in the following manner:
\begin{equation}
\begin{pmatrix} \label{bigmatrix}
\frac{D_1(x,t, \omega; \xi, \tau, \Omega)}{D(y,r,\theta ', \omega)} & \frac{D_2(x,t, \omega; \xi, \tau, \Omega)}{D(\bar{y},\bar{r},\bar{\theta}' \bar{\theta} '', \bar{\omega})}
\end{pmatrix}
\end{equation}
where the left submatrix corresponds to (\ref{flowout1}) and the right submatrix corresponds to (\ref{flowout2}).  We will first look at the concatenation of two submatrices in the above matrix, namely $\frac{D_1(\xi, \tau)}{D(\theta ')}$ and $\frac{D_2(\xi, \tau)}{D(\bar{\theta}', \bar{\theta}'')}$, which is

\begin{equation}
\begin{pmatrix}
\kappa^*d_{\theta '}(\nu - \sigma \omega) &  \kappa^*d_{\bar{\theta} '}(\bar{\nu} - \bar{\sigma} \bar{\omega}) &  \kappa^*d_{\bar{\theta} ''}(\bar{\nu} - \bar{\sigma} \bar{\omega})\\

 \kappa^*d_{\theta '}(- \sigma \omega) &  \kappa^*d_{\bar{\theta} '}(- \bar{\sigma} \bar{\omega}) &  \kappa^*d_{\bar{\theta} ''}( - \bar{\sigma} \bar{\omega})
\end{pmatrix}
\end{equation}
Doing a row reduction by subtracting a multiple of the $\tau$-row from the $\xi$-rows gives us 
\begin{equation}
\begin{pmatrix}
\kappa^*d_{\theta '}(\nu) &  \kappa^*d_{\bar{\theta} '}(\bar{\nu}) &  \kappa^*d_{\bar{\theta} ''}(\bar{\nu})\\

 \kappa^*d_{\theta '}(- \sigma \omega) &  \kappa^*d_{\bar{\theta} '}(- \bar{\sigma} \bar{\omega}) &  \kappa^*d_{\bar{\theta} ''}( - \bar{\sigma} \bar{\omega})
\end{pmatrix}
\end{equation}

As $\nu=\bar{\nu}$ on the intersection of our flowouts, the first $d_1$ columns are equal to the second $d_1$ columns of this matrix.  Moreover, as $\bar{\nu}$ forms a frame in the fibers of $N^*(S_{2+})$, the last $d_2$ columns are independent from the first $d_1$.  Hence, the rank of this concatenation of submatrices is $d_1+d_2$.  The rank of the differential in (\ref{flowout1}) is invariant meaning that if we do the same row reductions on the submatrix $\frac{D_1(\xi, \tau)}{D(\theta ')}$, then the rest of the rows in the differential $\frac{D_1(x,t, \omega; \xi, \tau, \Omega)}{D(y,r,\theta ', \omega)}$ will always give us rank $2n-d_1$.  Now, for every rank of $d_1$ given by this first submatrix, we get another $d_2$ for free by our observation.  Therefore (\ref{bigmatrix}) has rank at least $2n+d_2$ and $\Lambda_{1+}^a \cap \Lambda_{2+}^c$ is a clean intersection.

It is crucial for the intersection of the flowouts to be nonempty even after we avoid our bad set $\Sigma^1 \backslash \Sigma_2 \cup \Sigma^2 \backslash \Sigma_2$ for otherwise the standard microlocal analysis cannot be applied.  This can be checked by examining the parametrizations. 

\begin{rem}
The above calculations imply the intersection of the flowouts determine one of the initial surfaces, namely $S_2$.  
\end{rem}

\begin{rem}
While $\Lambda_{1+}^a$ and $\Lambda_{2+}^a$ are vector bundles over hypersurfaces, the clean intersection of these Lagrangians implies the underlying surfaces intersect tangentially.
\end{rem}

\subsection{Geometry associated to Radon transforms}

\noindent Let $R: \mathcal{E}'(\mathbb{R}^n) \rightarrow \mathcal{E}'(\mathbb{R} \times S^{n-1})$ be the Radon transform 
\begin{equation}
(Rf)(s, \theta,t, \omega) = \int_{x \cdot \theta=s} f(x,t,\omega) d\sigma(x)
\end{equation} 
where $d\sigma$ is normalized Lebesgue measure on the hyperplane $\{x \cdot \theta = s\}$.  $R$ is an elliptic Fourier integral operator \cite{Ho09}, with $R \in I^{-\frac{n-1}{2}}(C_R)$ where $C_R$ is the local canonical graph inside of $T^*(\mathbb{R} \times S^{n-1} \times \mathbb{R} \times S^{n-1}) \times T^*(\mathbb{R}^n \times \mathbb{R} \times S^{n-1})$ given by 
\begin{eqnarray}
\nonumber C_R &=& \{\big(x \cdot \theta, \theta, t, \omega; \sigma, \sigma i^*_{\theta}(x), \tau, \Omega \big),\big(x,t, \omega, \sigma \theta, \tau, \Omega \big): \\
\nonumber && (x,t, \omega) \in \mathbb{R}^{n+1} \times S^{n-1}, \theta \in S^{n-1}, \sigma \in \mathbb{R}-0, \tau \in \mathbb{R}, \\
&& \Omega \in T^*_{\omega}(S^{n-1})\}.
\end{eqnarray}  

An operator $F$ whose canonical relation is $C_R$ will be used when constructing the scattering kernel in Section \ref{sip}.  When we apply $C_R$ to $\Lambda_+$, $\Lambda_{1+}^a$, and $\Lambda_{2+}^a$, remembering that the application of canonical graphs to Lagrangians always satisfy the transverse intersection calculus, we get 
\begin{eqnarray}
\label{rad1} C_R \circ \Lambda_+ &=& \{(y \cdot \omega, \omega, y \cdot \omega, \omega; \sigma, - \sigma i^*_{\omega}(y), -\sigma, \sigma i^*_{\omega}(y)): y \in \mathbb{R}^n, \\
\nonumber && \omega \in S^{n-1}, \sigma \in \mathbb{R}-0\} 
\end{eqnarray}
\begin{eqnarray}
\nonumber &\cup& \{(-y \cdot \omega, -\omega, y \cdot \omega, \omega; \sigma, - \sigma i^*_{\omega}(y), \sigma, -\sigma i^*_{\omega}(y)): y \in \mathbb{R}^n, \\
\nonumber && \omega \in S^{n-1},  \sigma \in \mathbb{R}-0\},
\end{eqnarray}
\begin{eqnarray}
\label{rad2} C_R \circ \Lambda_{1+}^a &=& \{(r- y \cdot \frac{(\nu - \sigma \omega)}{\sigma}, -\frac{(\nu - \sigma \omega)}{\sigma}, y \cdot \omega + r, \omega; \sigma, \\
\nonumber &&\sigma i^*_{-\frac{(\nu - \sigma \omega)}{\sigma}}(y - r(\nu -\sigma \omega)), \sigma, -\sigma i^*_{\omega}(y)): \\
\nonumber && (y, \nu) \in N^*(S_{1}), \omega \in S^{n-1}, r \in \mathbb{R}, \sigma = \frac{|\nu|^2|}{2 \nu \cdot \omega} \text{ with } \nu \cdot \omega \neq 0\} 
\end{eqnarray}
\begin{eqnarray}
\nonumber &\cup& \{(-r + y \cdot \frac{(\nu - \sigma \omega)}{\sigma}, \nu - \sigma \omega, y \cdot \omega + r, \omega; \sigma, \\
\nonumber && \sigma i^*_{-\frac{(\nu - \sigma \omega)}{\sigma}}(y - r(\nu -\sigma \omega)), \sigma, -\sigma i^*_{\omega}(y)):\\
\nonumber && (y, \nu) \in N^*(S_{1}), \omega \in S^{n-1}, r \in \mathbb{R}, \sigma = \frac{|\nu|^2|}{2 \nu \cdot \omega} \text{ with } \nu \cdot \omega \neq 0\},
\end{eqnarray}
\begin{eqnarray}
\label{rad3} C_R \circ \Lambda_{2+}^c &=& \{(r- y \cdot \frac{(\nu - \sigma \omega)}{\sigma}, -\frac{(\nu - \sigma \omega)}{\sigma}, y \cdot \omega + r, \omega; \sigma,\\
\nonumber && \sigma i^*_{-\frac{(\nu - \sigma \omega)}{\sigma}}(y - r(\nu -\sigma \omega)), \sigma, -\sigma i^*_{\omega}(y)): \\
\nonumber&& (y, \nu) \in N^*(S_{2}), \omega \in S^{n-1}, r \in \mathbb{R}, \sigma = \frac{|\nu|^2|}{2 \nu \cdot \omega} \text{ with } \nu \cdot \omega \neq 0\}
\end{eqnarray}
\begin{eqnarray*}
\nonumber &\cup& \{(-r + y \cdot \frac{(\nu - \sigma \omega)}{\sigma}, \nu - \sigma \omega, y \cdot \omega + r, \omega; \sigma, \\
\nonumber && \sigma i^*_{-\frac{(\nu - \sigma \omega)}{\sigma}}(y - r(\nu -\sigma \omega)), \sigma, -\sigma i^*_{\omega}(y)):\\
\nonumber && (y, \nu) \in N^*(S_{2}), \omega \in S^{n-1}, r \in \mathbb{R}, \sigma = \frac{|\nu|^2|}{2 \nu \cdot \omega} \text{ with } \nu \cdot \omega \neq 0\}.
\end{eqnarray*}
We note that as $C_R$ is a local canonical graph, i.e. locally the graph of a symplectomorphism, $(C_R \circ \Lambda_{1+}^a, C_R \circ \Lambda_{2+}^c)$ is a pair of Lagrangians intersecting cleanly in codimension $d_2$.  

\subsection{Geometry associated to pullbacks}

\noindent For $t_0 >> 0$, the mapping $\rho_{t_0}: \mathbb{R} \times S^{n-1} \times S^{n-1} \rightarrow \mathbb{R} \times S^{n-1} \times \mathbb{R} \times S^{n-1}$, given by $\rho_{t_0}=(s, \theta, \omega) = (t_0 + s, \theta, t_0, \omega)$, induces the restriction mapping
\begin{equation}
\rho^*_{t_0}: \mathcal{D}_{\rho}'(\mathbb{R} \times S^{n-1} \times \mathbb{R} \times S^{n-1}) \rightarrow \mathcal{D}'( \mathbb{R} \times S^{n-1} \times S^{n-1}),
\end{equation}
where $\mathcal{D}_{\rho_{t_0}}'$ is the space of distributions whose wavefront is disjoint from the normals of $\rho_{t_0}$.  We use this space in order to make this restriction map well-defined.  We drop the subscript $t_0$ from here on, as it is understood what $\rho_{t_0}$ means.  It follows that $\rho^*$ is Fourier integral operator and $\rho^* \in I^{\frac{1}{4}}(C_{\rho})$, where 
\begin{eqnarray}
\nonumber C_{\rho} &=& \{(s, \phi, \omega; \tau, \Phi, \Omega),(t_0 + s, \phi, t_0, \omega; \tau, \Phi, \eta, \Omega): s \in \mathbb{R}, \phi \text{ and } \omega \in S^{n-1}, \\
&& (t_0 + s, \phi, t_0, \omega; \tau, \Phi, \eta, \Omega) \in T^*(\mathbb{R} \times S^{n-1} \times \mathbb{R} \times S^{n-1}) \backslash 0\}
\end{eqnarray}
It was proven in \cite{GU93} that the compositions of $C_{\rho}$ with the Lagrangians in (\ref{rad1}),(\ref{rad2}), and (\ref{rad3}) are all transversal.  In fact, the second components in (\ref{rad1}),(\ref{rad2}), and (\ref{rad3}) vanish after the application of $C_{\rho}$.
\begin{defn}
The \underline{peak Lagrangian} is defined as
\begin{eqnarray}
\hat{\Lambda}_+ &=& C_{\rho} \circ C_{R} \circ \Lambda_+ \\
\nonumber &=& \{(y, \omega, \omega; \sigma, -\sigma i^*_{\omega}(y), \sigma i^*_{\omega}(y)): y \in \mathbb{R}^n, \omega \in S^{n-1}, \sigma \in \mathbb{R} \backslash 0\}.
\end{eqnarray}
\end{defn}

\begin{defn}
The \underline{reflected Lagrangians} are defined as 
\begin{eqnarray}
\label{ref1} \hat{\Lambda}_{1+}^a &=& C_{\rho} \circ C_{R} \circ \Lambda^a_{1+} = \{(-y \cdot (\frac{(\nu - \sigma \omega)}{\sigma} + \omega),- \frac{(\nu - \sigma \omega)}{\sigma}, \omega; \\
\nonumber && \sigma, -\sigma i^*_{\frac{(\nu - \sigma \omega)}{\sigma}}(y), \sigma i^*_{\omega}(y)): (y, \nu) \in N^*(S_1), \omega \in S^{n-1}, \\
\nonumber && r \in \mathbb{R}, \sigma = \frac{|\nu|^2|}{2 \nu \cdot \omega} \text{ with } \nu \cdot \omega \neq 0\}
\end{eqnarray}
and
\begin{eqnarray}
\label{ref2} \hat{\Lambda}_{2+}^c &=& C_{\rho} \circ C_{R} \circ \Lambda^c_{2+} = \{(-y \cdot (\frac{(\nu - \sigma \omega)}{\sigma} + \omega),- \frac{(\nu - \sigma \omega)}{\sigma}, \omega; \\
\nonumber && \sigma, -\sigma i^*_{\frac{(\nu - \sigma \omega)}{\sigma}}(y), \sigma i^*_{\omega}(y)): (y, \nu) \in N^*(S_2), \omega \in S^{n-1}, \\
\nonumber && r \in \mathbb{R}, \sigma = \frac{|\nu|^2|}{2 \nu \cdot \omega} \text{ with } \nu \cdot \omega \neq 0\}
\end{eqnarray}
\end{defn}
\noindent A direct calculation shows $\hat{\Lambda}_{1+}^a \cap \hat{\Lambda}_{2+}^a$ is a smooth submanifold of codimension $d_2$ inside each of the reflected Lagrangians. 

The proof that $\hat{\Lambda}_{1+}^a \cap \hat{\Lambda}_{2+}^a$ is clean is almost identical to that which shows $\Lambda_{1+}^a \cap \Lambda_{2+}^c$ is clean.  The concatenation of the differentials of the parametrizations in (\ref{ref1}) and (\ref{ref2}) will have rank $2n+d_2$ if the concatenation of the submatrices $\frac{D_1(\phi, \tau)}{D(\theta ')}$ and $\frac{D_2(\phi, \tau)}{D(\bar{\theta '},\bar{\theta ''})}$ has rank $d_1+d_2$.  This follows immediately after realizing that $\frac{D_1(\phi, \tau)}{D(\theta ')}$ is just the composition of $\frac{D_1(\xi)}{D(\theta ')}$ from (\ref{lambda1+}), which has rank $d_1$, with the polar coordinate map, which has an invertiable differential away from the origin.  The same holds true for $\frac{D_2(\phi, \tau)}{D(\bar{\theta} ',\bar{\theta} '')}$ in relation to $\frac{D_2(\xi)}{D(\bar{\theta} ', \bar{\theta} '')}$ from (\ref{lambda2+}), which has rank $d_1+d_2$.  Moreover, $\frac{D(\xi)}{D(\theta ')}$ has as its columns the first $d_1$ columns of $\frac{D(\xi)}{D(\bar{\theta} ', \bar{\theta} '')}$.  The remaining reasoning is similar.

\begin{rem}
Once again, the clean intersection calculations indicate that the reflected Lagrangians should give us substantial information about $S_1$ and $S_2$.  This will be justified in the final section of this paper.
\end{rem}

\section{Multiplication by a nested conormal} \label{mnc}

\noindent The operator $M_q$, which is multiplication by $q \in I^{M_1, M_2}(S_1, S_2)$, has wavefront set in the three canonical relations
\begin{eqnarray}
\triangle_{T^*(\mathbb{R}^n \times \mathbb{R} \times S^{n-1}) \times T^*(\mathbb{R}^n \times \mathbb{R} \times S^{n-1})},
\end{eqnarray}
\begin{eqnarray} 
 C_{S_1} & = & \{ (x,t,\omega ; \xi, \tau, \Omega), (x,t,\omega ; \xi + D_x((h_i)_{i=1}^{d_1})\cdot \theta, \tau, \Omega): \\
\nonumber & & (x,t,\omega ; \xi, \tau, \Omega) \in T^*(\mathbb{R}^n \times \mathbb{R} \times S^{n-1})_{|S_i}, \text{ and } \theta \in \mathbb{R}^{d_1} \backslash 0\},
\end{eqnarray}
 and 
\begin{eqnarray} 
C_{S_2} & = & \{ (x,t,\omega ; \xi, \tau, \Omega), (x,t,\omega ; \xi + D_x((h_i)_{i=1}^{d_1+d_2})\cdot \theta, \tau, \Omega): \\
\nonumber & & (x,t,\omega ; \xi, \tau, \Omega) \in T^*(\mathbb{R}^n \times \mathbb{R} \times S^{n-1})_{|S_2}, \text{ and } \theta \in \mathbb{R}^{d_1+d_2} \backslash 0\}.
\end{eqnarray}

An operator of this type no longer falls into the class of operators whose Schwartz kernels are paired Lagrangians; in fact, the above canonical relations form a triple of pairwise cleanly intersecting Lagrangians.  The lack of a developed theory for the compositions of such operators requires us to take some different approaches when analyzing its mapping properties on distributions.

$M_q$ is technically defined on the space of distributions $\mathcal{D}$ that has the property $(C_{S_1} \circ WF(v)) \cap 0 = \emptyset$ and $(C_{S_2} \circ WF(v)) \cap 0 = \emptyset$, for all $v \in \mathcal{D}$.  However, the multiplicative results and Sobolev mapping properties of this section circumvent this technicality. 

\subsection{Action on spaces of distributions}

\noindent The following is a result from \cite{GU93}.

\begin{lem} \label{singlemult}
Let $X$ be a manifold and $Y, Z \subset X$ be submanifolds with $Y \pitchfork Z$, then 
\begin{equation} \label{trans mult}
I^{\mu}(Y) \cdot I^{\mu'}(Z) \subset I^{\mu, \mu'}(Y, Y \cap Z) \oplus I^{\mu, \mu'}(Z, Y \cap Z).
\end{equation}
Moreover, if $u \in I^{\mu}(Y)$ satisfies supp \hskip.1cm $u$ $\subset Y$, then 
\begin{equation}
u \cdot I^{\mu}(Z) \subset I^{\mu, \mu '}(Y, Y \cap Z).
\end{equation}
\end{lem}

\noindent We need another lemma stating the multiplication properties of distributions associated to nested pairs of submanifolds.

\begin{lem} \label{doublemult}
Suppose $Y_+, Y_1,$ and $Y_2$ are submanifolds of $X$ such that $Y_2 \subset Y_1$, $Y_+ \pitchfork Y_1$, and $Y_+ \pitchfork Y_2$.  If $u_1 \in \in I^{\mu}(Y_+), u_2 \in  I^{M_1, M_2}(Y_1, Y_2)$ with supp $u_1$ $\subset Y_+$ and $u_2$ supported microlocally near $N^*(Y_1) \cap N^*(Y_2)$, then 
\begin{equation}
u_1u_2 \in I^{\mu, M_1, M_2}(Y_+,Y_1 \cap Y_+ ,Y_2 \cap Y_+).
\end{equation}
\end{lem}

\noindent Note that this result is a direct analog of the previous lemma for nested conormal distributions.

\noindent \underline{Proof}: We can introduce local coordinates $(x',x'',x''') \in \mathbb{R}^{d_1} \times^{d_2} \mathbb{R} \times \mathbb{R}^{n-d_1-d_2}$ near the point $x_0 \in Y_+ \cap Y_2$  such that $x_0=0$, $Y_+ = \{x'=0\}$, $Y_1 = \{x''=0\}$, and $Y_2 = \{x''=x'''=0\}$.  Under these assumptions, 
\begin{equation*}
u_1(x) = \int_{\mathbb{R}^{d_1}} e^{i(x' \cdot \xi ')} a(x; \xi ') d \xi '
\end{equation*}
for $a \in S^{\mu}(X \times (\mathbb{R}^{d_1} \backslash 0))$, and $u_2$ as the representation 
\begin{equation*}
u_2(x) = \int_{\mathbb{R}^{d_1+d_2}} e^{i(x'' \cdot \xi '' + x''' \cdot \xi ''')} b(x; \xi '', \xi ''') d \xi '' d \xi '''
\end{equation*}
for $b \in S^{M_1, M_2}(X \times (\mathbb{R}^{d_1} \backslash 0) \times \mathbb{R}^{d_2})$.  Multiplying the distributions leads us to the oscillatory integral
\begin{equation}
(u_1u_2)(x) = \int_{\mathbb{R}^{d_1+d_2}} e^{i(x' \cdot \xi ' + x'' \cdot \xi '' + x''' \cdot \xi ''')} a(x; \xi ') b(x; \xi '', \xi ''') d \xi ' d \xi '' d \xi '''.
\end{equation}
since $u_2$ being supported near the aforementioned intersection is equivalent to having $b$ localized to the region $= \{\xi '' | \gtrsim |\xi '''| \}$. There are three regions in which $\xi '$ can lie: I $= \{ |\xi '| \gtrsim | \xi '' | \gtrsim |\xi '''| \}$, II $= \{ |\xi ''| \gtrsim | \xi ' | \gtrsim |\xi '''| \}$, and III $= \{ |\xi ''| \gtrsim | \xi ''' | \gtrsim |\xi '| \}$.  Computing the size of the Japanese brackets along with their orders, it follows that regions II and III correspond to the nested conormal $I(Y_1 ,Y_2 \cap Y_+)$ and the nested triple conormal $I(Y_1,Y_1 \cap Y_+ ,Y_2 \cap Y_+)$, respectively and of certain orders.  A general fact involving the multiplication of distributions from \cite{Ho71} says that 
\begin{equation*}
WF(u_1u_2) \subset N^*(Y_+) \cup N^*(Y_+ \cap Y_1) \cup N^*(Y_+ \cap Y_2)
\end{equation*}
as supp $u_1$ $\subset Y_+$.  All three terms above lie inside the class $I^{\mu, M_1, M_2}(Y_+,Y_1 \cap Y_+ ,Y_2 \cap Y_+)$, defined in Section \ref{ntc}.  Moreover, the region in which this integral has the singularity corresponding to the triple intersection is I because of $ab$ satisfying a triple symbol-valued symbol estimate, namely $ab \in S^{\mu, M_1, M_2}$.  $\blacksquare$

\subsection{Sobolev estimates}

\noindent In this section, we prove $L^2$ estimates for the operator that multiplies by a nested conormal distribution.  Previous estimates for multiplication by a standard conormal distribution were obtained in \cite{GU90} to describe operators with two canonical operators, one being a flowout relation and the other being the diagonal relation.  This work relied heavily on a composition calculus for certain kinds of paired Lagrangians developed by Antoniano and Uhlmann \cite{AntUhl85} and follows Hormander's method for obtaining $L^2$ estimates for Fourier integral operators associated to canonical graphs.  

As stated in Section \ref{lscr}, $M_q$ is now associated to three canonical relations, two of which are  flowout relations.  The lack of a formal composition calculus for operators associated to three intersecting Lagrangians forces us to take another approach.  We will make use of an observation of Melrose \cite{Melnotes89} (see also \cite{GU90}) that elements of $I^{M_1, M_2}(Y_1, Y_2)$ can be decomposed into a sum of two conormal distributions with $(\frac{1}{2}, \frac{1}{2})$ symbols; this is called a \underline{parabolic decomposition}.  From here, we use an orthogonality argument to essentially retrieve the same estimate as for multiplication by a standard conormal distribution, but in the $(\frac{1}{2}, \frac{1}{2})$ case.  For more on $(\frac{1}{2}, \frac{1}{2})$ symbols, see \cite[Chapter 7]{Ste93}.

\begin{prop} \label{paradecomp} \cite{Melnotes89}
Let $Y_1 \subset Y_2 \subset X$ with dim $X = n$, codim $Y_1 = d_1$ and $codim Y_2 = d_1 + d_2$.  Then for $-d_2 < M_2$,
\begin{equation}
I^{M_1, M_2}(Y_1, Y_2) \subset I^{M_1 + \frac{M_2}{2} + \frac{d_2}{2}}_{\frac{1}{2}, \frac{1}{2}}(Y_1) +  I^{M_1 + \frac{M_2}{2}}_{\frac{1}{2}, \frac{1}{2}}(Y_2),
\end{equation}
while for $M_2 \leq - d_2$
\begin{equation}
I^{M_1, M_2}(Y_1, Y_2) \subset I^{M_1}_{\frac{1}{2}, \frac{1}{2}}(Y_1) +  I^{M_1 + \frac{M_2}{2}}_{\frac{1}{2}, \frac{1}{2}}(Y_2),
\end{equation}
where we continue to use conormal notation.  Here, the spaces on the left hand side are conormal distributions with symbols of type $(\frac{1}{2}, \frac{1}{2})$. 
\end{prop}

\noindent We recall a version of the Cotlar-Knapp-Stein lemma \cite{Ste93} that will be used in proving our Sobolev estimates.

\begin{lem} \label{cotlar}
Let $j \in \mathbb{Z}^r$, $T = \sum_{j \in \mathbb{Z}^r} T_j$ with $\{T_j\}_{j \in \mathbb{Z}^r}$ bounded sequence of operators on $L^2$, and $\{ \gamma(j) \}_{j \in \mathbb{Z}^r}$ be a multiparameter sequence of positive numbers such that $A = \sum_{j \in \mathbb{Z}^r} \gamma(j) < \infty$.  If \begin{eqnarray*}
\|T_i^*T_j\| \leq (\gamma(i-j))^2 \\
\|T_iT^*_j\| \leq (\gamma(i-j))^2
\end{eqnarray*}
then 
\begin{equation*}
\|T\| \leq A.
\end{equation*}
\end{lem}

\begin{thm} \label{sobomult}
Let $M_q$ be the operator that multiplies by $q \in  I^{M}_{\frac{1}{2}, \frac{1}{2}}(S)$, where $S$ has codim$=d$ in $X$. For $M=-d+\alpha$ with $0 \leq \alpha < d$, $M_q$ maps $H^s$ to $H^{s-\alpha}$.
\end{thm}

\begin{rem}
The assumption that $\alpha < d$ helps avoid multiplication by a distribution that has as strong of a singularity as the delta function.
\end{rem}

\noindent \underline{Proof:}  As Sobolev spaces are diffeomorphism invariant, we can apply an elliptic Fourier integral operator associated to a canonical transformation that turns $q$ into 

\begin{equation}
\int e^{ix' \cdot \theta ' }a(x;\theta ') d \theta '
\end{equation}
where $a \in S^M_{\frac{1}{2}, \frac{1}{2}}$ \cite{Ho71}.  

Notice the kernel of $M_q$ is equal to $q(y) \delta(x-y)$.  An interesting feature of this operator is that the kernel is also equal to $q(x) \delta(x-y)$.  Also note $M_q^*$ corresponds to multiplication by $\bar{q}$.  We take all these facts into account when composing the two operators.

Let $\{\chi_j : j \geq 0 \}$ be a non-homogeneous dyadic partition of unity on $\mathbb{R}$ such that
\begin{eqnarray*}
&&\text{ supp }(\chi_0) \subset \{|t| \leq 2\} \\
&&\text{ supp }(\chi_j) \subset \{2^{j-1} \leq |t| \leq 2^{j+1}\}, j \geq 1, \text{ and } \\
&&|\chi_j^{m}| \leq C_m2^{-mj}, \text{ for all } m \geq 1.
\end{eqnarray*}

In the language of Lemma \ref{cotlar}, the terms $M_{q_{j,k}}$ will have Schwarz kernels 
\begin{equation} 
\int e^{i[(x-z) \cdot \xi + z' \cdot \theta]} \chi_k(| \xi |) a_j(z;\theta) d \theta d \xi .
\end{equation}
Set $\chi_{j,j',k,k'}=\chi_j(|\theta|) \chi)_{j'}(|\tau|) \chi_k(|\xi|) \chi_{k'}(|\eta|).$  Since the $q_{j,k}$'s are $C^{\infty}$ functions, the compositions $M_{q_{j,k}}M_{q_{j',k'}}^*$ are well-defined and equal to $M_{q_{j,k}}^*M_{q_{j',k'}}$ as the multiplication operators are normal.  The Schwarz kernels for these operators are

\begin{eqnarray} \label{normaloperator}
K_{j,j',k,k'}(x,y) = \int && e^{i[(x-z) \cdot \xi + z' \cdot \theta - (y-z) \cdot \eta - z' \cdot \tau]} \times \\
\nonumber && \chi_k(|\xi|) \chi_{k'}(|\eta|) a_j(z;\theta)a_{j'}(z;\tau) dz d \theta d \tau d \xi d \eta
\end{eqnarray}
where $a_j = \chi_j \cdot a$ and $a_{j'}= \chi_{j'} \cdot a$.  The numbers $\gamma(j)$ will be of size $2^{-j}$ as a result of our partition of unity.  This kind of decomposition of an operator is commonly referred to as \underline{Littlewood-Paley decomposition}; see \cite{Ste93}. For brevity in our upcoming calculations, we refer to the phase in (\ref{normaloperator}) as $\Phi$.

Consider the operator $L_z = I - \triangle_z$, where $\triangle_z$ denotes the Laplacian in the $z$-variable.  It follows that $L^t_z=L$.  Notice $$\frac{(I - \triangle_z)^{N_z}}{(\langle \xi + (\theta,0) - \eta - (\tau,0) \rangle)^{N_z}}e^{i\Phi}=e^{i\Phi},$$ with the expression on the bottom being $\langle \nabla_z \Phi \rangle^{N_z}$.  It is the size of this Japanese bracket that will determine how we get our estimate in each region of integration.  There are several cases we must consider:

\vskip.2cm \noindent \underline{Case 1:} Single elliptic variable.  

\noindent This region of integration corresponds to the part of the operator that is microlocally supported far away from the diagonal.  Without loss of generality, because of symmetry between $\xi$ and $\tau$ (likewise for $\theta$ and $\tau$), we can let the elliptic variable be $\xi$.  In this region, $\langle \nabla_z\Phi \rangle \approx \langle \xi \rangle$.  Remembering that a spatial differentiation of symbols of type $(\frac{1}{2},\frac{1}{2})$ loses $\langle \xi \rangle  ^{\frac{1}{2}}$, integration by parts $N_z$ times with our differential operator $L_z$ gives 

\begin{align}
\nonumber \int e^{i\Phi} \langle \nabla_z \Phi \rangle ^{-N_z}[ (L_z^T)^{N_z}\chi_k(|\xi|) \chi_{k'}(|\eta|) a(z;\theta)a(z;\tau)] dz d \theta d \tau d \xi d \eta \\
\nonumber \lesssim \int \chi_k(|\xi|) \chi_{k'}(|\eta|) \chi_j(|\theta|) \chi_{j'}(|\tau|) <\xi>^{-2N_z}<\xi>^{\frac{N_z}{2}}<\xi>^{\frac{N_z}{2}} dz d \theta d \tau d \xi d \eta \\
\lesssim 2^{M(j+j')+d(j+j')+n(k+k) - N_zk} \times Vol(\text{support of }z).
\end{align}
The above estimate is still rough as we have not utilized the integration in $z$ and the later integrations in $x$ or $y$ demanded by the application of Young's inequality; see \cite{Foll95}.

If we let $L_{\xi}$ and $L_{\eta}$ be defined analogously to $L_z$, with the differentiations happening in the variables determined by the subscripts, then integration by parts $N_{\xi}$ and $N_{\eta}$ times in $\xi$ and $\eta$, respectively, gives us 
\begin{eqnarray*}
\int e^{i\Phi} &\times& \langle 2^k(x-z) \rangle ^{-N_{\xi}} (L_{\xi}^T)^{N_{\xi}} \langle 2^{k'}(y-z) \rangle ^{-N_{\eta}} (L_{\eta}^T)^{N_{\eta}} \times  \\
&& \langle \nabla_z \Phi \rangle ^{-N_z} (L_z^T)^{N_z}[\chi_k(|\xi|) \chi_{k'}(|\eta|) a_j(z;\theta)a_{j'}(z;\tau)] dz d \theta d \tau d \xi d \eta
\end{eqnarray*}
Remembering $(L_{\xi}^T)^{N_{\xi}} (L_{\xi}^T)^{N_{\eta}} \langle \xi + (\theta,0) - \eta - (\tau,0) \rangle^{-N_z} \lesssim \langle \xi + (\theta,0) - \eta - (\tau,0) \rangle^{-N_z}$ we get a better size estimate of 
\begin{align} \label{fin size est}
I_1 \times 2^{M(j+j')+d(j+j')+n(k+k) - (2N_z)k}, 
\end{align}
where $I_1 = \int \langle 2^k(x-z) \rangle ^{-2N_{\xi}}  \langle 2^{k'}(y-z) \rangle ^{-2N_{\eta}} dz$.  Since $\frac{1}{2}|x-y| \leq |x-z|$ or $\frac{1}{2}|x-y| \leq |y-z|$, it follows 
\begin{eqnarray}
\nonumber I_1 \lesssim I_2 &=& \int \langle 2^{k-1}(x-y) \rangle ^{-2N_{\xi}}  \langle 2^{k'}(y-z) \rangle ^{-2N_{\eta}} dz + 
\\ && \int \langle 2^k(x-z) \rangle ^{-2N_{\xi}}  \langle 2^{k'-1}(x-y) \rangle ^{-2N_{\eta}} dz.
\end{eqnarray}
Using this inequality and completing the $z$ integration in $I_2$ followed by an integration in $x$ or $y$, we see 
\begin{equation}
\|M_{j,j',k,k'}M^*_{j,j',k,k'}\|_2 \lesssim 2^{M(j+j')+d(j+j')-(2N_z)k}.
\end{equation}
Letting $N_z > M + d$, we get rapid decay of the norm.

This kind of behavior is typical in regions far away from the diagonal of phase space.  The following cases will show that such freedom for $M$ is not allowed when phase variables are closer to the diagonal.

\vskip.2cm \noindent \underline{Case 2}: $\{|\xi| \gtrsim |\eta|, |\theta| \approx |\tau|\}$ or $\{|\xi| \approx |\eta|, |\theta| \gtrsim |\tau|\}$

\noindent Suppose we are in the first region of integration, $\langle \nabla_z \Phi \rangle \approx \langle \xi - \eta \rangle \approx \langle \xi \rangle$.  The situation when $\xi$ is elliptic to $\theta$ or $\tau$ is covered in Case 1.  Otherwise, the gain obtained by $\langle \nabla_z \Phi \rangle$ does not overwhelm the loss of $\langle \theta \rangle ^{\frac{N_z}{2}}\langle \tau \rangle ^{\frac{N_z}{2}}$.  The reasoning for the second region of integration is identical.  Hence, the upper bound for $K_{j,j',k,k'}(x,y)$ becomes
\begin{align}
\label{trivial} I_1 \times 2^{M(j+j')+d(j+j')+n(k+k')}
\end{align}

Our only option left is to use the pseudodifferential operator part of $M_q$ to get a gain in $k$ (or $k'$) after integration in $z$ followed by an integration in $x$ or $y$ to get the additional gain in $k'$ (or $k$) and finish the application of Young's inequality, just as in Case 1.  This leaves us setting $M < -d$ in order to sum along these indices.

Because of the aforementioned symmetry in Case 1, the analysis in this region covers the situation when the elliptic variable in a pair is switched.

\vskip.2cm \noindent \underline{Case 3}:  $\{|\xi| \approx |\eta|, |\theta| \approx |\tau|\}$

\noindent The quantity $\langle \nabla_z \Phi \rangle$ is bounded above by 1 as
\begin{align*}
\langle \nabla_z \Phi \rangle = \langle \xi + (\theta,0) - \eta - (\tau,0) \rangle \approx \langle \xi  - \eta \rangle \approx 1,
\end{align*}
implying that the only gain we get comes from the pseudodifferential operator part of $M_q$, just like in Case 2.  Following those arguments again, we must assume that $M < -d$.

\vskip.2cm \noindent \underline{Case 4}: $\{|\xi| \approx |\eta| \approx |\theta| \approx |\tau|\}$.

\noindent We are now on the diagonal and must estimate the operators $M_{j,j}$ whose kernels are of the form 
\begin{align*}
K_{j,j}(x,y)=\int e^{i[x' \cdot \theta + (x-y)\cdot \xi]}\chi_j(\xi) a_j(x;\theta) d \theta d \xi.
\end{align*}
Since there are no spatial integrations of which to take advantage, we can only integrate by parts $N_{\xi}$ times in $\xi$, giving us the estimate
\begin{align*}
\int \langle 2^j(x-y) \rangle^{-N_{\xi}} (L_{\xi}^T)^{N_{\xi}}[\chi_j(|\xi|)a_j(x;\theta)] d \theta d \xi \lesssim \langle 2^j(x-y) \rangle^{-N_{\xi}} 2^{Mj + dj + nj}. 
\end{align*}
Integrating in $x$ or $y$ gives an additional gain of $2^{-nj}$, forcing us again to have $M < -d$ for summability.

\vskip.2cm \noindent \underline{Case 5}: All other regions.

\noindent It worth noting that outside of Case 1, all of the analysis done was dependent upon whether the gain obtained from the spatial gradient $\nabla_z \Phi$ overwhelmed the loss from the symbol.  In each case, integration by parts in the phase variables of the pseudodifferential part of the kernel and the last spatial integration in $x$ or $y$ always left us with 
\begin{equation}
|K_{j,j',k,k'}(x,y)| \lesssim I_1 \times 2^{M(j+j')+d(j+j')+n(k+k')},
\end{equation}
before the $z$ integration.  Hence, all of the above analysis, including in Case 1, could have been done this way.  The effort was made in Case 1 to show that microlocalization far away from the diagonal gives rapid decay, as expected.

Therefore, by following the steps in Case 2 or after, any other region of integration requires us to have $M < -d$ like in the previous cases.  If $M = -d + \alpha$, we see the norms grow at a rate of $2^{j \alpha}$.  For a fixed $\varepsilon > 0$, dividing out by $2^{j \alpha + \varepsilon}$ gives the $H^{-\alpha-\varepsilon}$ norm of $M_{q_{j,k}}M_{q_{j',k'}}^*$. $\blacksquare$

\begin{cor}
As in Proposition \ref{paradecomp} with $\alpha < d_1$, we have $M_q: H^s(X) \rightarrow H^{s - \alpha}$ for 

\begin{center}
$\begin{cases}
M_1 + \frac{M_2}{2} < -d_1 - d_2 + \alpha, & \text { for } M_2 > -d_2 \text{ and } M_1 < -d_1 - \frac{d_2}{2} + \alpha\\
M_1 + \frac{M_2}{2} < -d_1 - d_2 + \alpha, & \text{ for } M_2 \leq -d_2 \text{ and } M_1 < -d_1 + \alpha
\end{cases}$
\end{center}
\end{cor}

\noindent \underline{Proof}: Applying the parabolic cutoff to $M_q$, we get the sum $M_{q_1} + M_{q_2}$ with $q_1$ and $q_2$ being of type $(\frac{1}{2}, \frac{1}{2}).$  We apply our last proposition to obtain the necessary restriction on our orders.  $\blacksquare$

\begin{rem}
Notice that for $M_2$ sufficiently negative, $M_1$ is allowed to be greater than $-d_1$.  By estimates on the growth rate of nested conormal distributions in Section \ref{pre}, this regime of orders allows $q$ to blow up, i.e. be unbounded near $S_1$.  Moreover, using the fact that
\begin{equation*}
\bigcap_{M_2 \in \mathbb{R}} I^{M_1,M_2}(Y_1, Y_2) = I^{M_1}(Y_1),
\end{equation*}
 we have generalized an estimate used in \cite{GU90} in the case of multiplication operators.  This theorem will be stated in the next section.
\end{rem}

\section{Parametrix to the wave equation} \label{pwo}

\noindent Parametrices for operators $P(x,D)$ of real principal type, such as $\square = \partial^2_t - \triangle$, were originally considered by Duistermaat and H\"ormander \cite{DuisHo72} and were shown by Melrose and Uhlmann \cite{MelUhl79} to be operators whose Schwartz kernels are paired Lagrangians, whose wave fronts are contained inside the union of the diagonal relation and the flowout from the diagonal under the Hamiltonian vector field of the principal symbol of $P(x,D)$.  The fundamental solution $\square^{-1}$, which is an exact parametrix, lies in the space $I^{-\frac{3}{2}, -\frac{1}{2}}(\triangle, C_{\square})$, where 
\begin{eqnarray} \label{flowout}
\nonumber C_{\square} &=& \{\big(x,t, \omega; \xi, |\xi|, \Omega\big),\big(x + s \frac{\xi}{|\xi|}, s, \omega ; \xi, |\xi|, \Omega \big) :\\ 
&& (x,t, \omega) \in \mathbb{R}^n \times \mathbb{R} \times S^{n-1}, \xi \in \mathbb{R}^n \backslash 0, \Omega \in T^*_{\omega}(S^{n-1}) \}
\end{eqnarray} 
is the canonical relation obtained by flowing out from the light cone $\{ |\xi| = |\tau| \}$ along the vector field $H_{\square}= \xi \cdot \frac{\partial}{\partial x} + \tau \frac{\partial}{\partial t}$.

The Schwartz kernel of $\square^{-1}$ can be written in the form 
\begin{equation} \label{boxinverse}
\int e^{i[(x-y)\cdot \xi + (t-s)|\xi| + (t-s) \rho]} a(x,t,y,s; \xi, \rho) d \xi d \rho
\end{equation} 
where the function appearing in the exponential is a multiphase parametrizing the associated canonical relations and $a \in S^{-\frac{1}{2},-\frac{3}{2}}(\mathbb{R}^{n+1} \times \mathbb{R}^{n+1} \times (\mathbb{R}^{n} \times \mathbb{R}) \backslash 0)$.  This generalized Fourier integral operator appeared in \cite{DuisHo72} but was not systematically treated until the work \cite{MelUhl79} and Guillemin and Uhlmann \cite{GuiUhl81}.  Due to the presence of the diagonal relation, which is the wavefront set of standard pseudodifferential operators, operators of this type are given the name of ``pseudodifferential operators with singular symbols"; an in-depth study on operators of this type is done in \cite{GU90}.  We will recall some mapping properties of such operators, both on Lagrangian-type distributions and Sobolev spaces, from this article which will be of use when solving our inverse problem.

The study of operators similar to $\square^{-1}$ from a microlocal perspective has been carried out in various geometrical settings and can more recently be found, for example, in the work of Baskin \cite{Bask10} on scattering theory on de Sitter space and anti-de Sitter space.  In \cite{Josh98a}, Joshi develops a symbolic calculus similar to that of Melrose and Uhlmann to construct complex powers of $\square$ on Riemannian manifolds and develop parametrices for the resulting operators. 

\subsection{Action on spaces of distributions}
~
We now consider the mapping properties of a parametrix for a pseudodifferential operator of real principal type on various kinds of Lagrangian distributions.  Of course, our intended application is for the d'Alembertian on $\mathbb{R}^{n-1} \times S^{n-1}$, but because of the complexity of certain phase functions in our coordinate system, we prove the results in more generality in order to utilize the more computationally convenient normal forms.  

Let $P(x,D)$ be an $m$th order pseudodifferential operator with real homogeneous principal symbol of classical type $p_m(x, \xi)$.  If $C_P = \{ (x; \xi),(y ; \eta): (x; \xi) \in char(P), (y, \eta) \in \Xi_{(x; \xi} \}$ where $\Xi_{(x; \xi)}$ is a bicharacteristic of $P(x,D)$, it follows from \cite{MelUhl79} that the parametrix $Q$ for $P(x,D)$ lies in the class $I^{\frac{1}{2}-m, - \frac{1}{2}}( \triangle, C_{P})$.  The following results are from \cite{GU93}.

\begin{prop} \label{singleflow}
Suppose $\Lambda_0 \subset T^*(X) \backslash 0$ is a conic Lagrangian intersecting $char(P)$ transversally and such that each bicharacteristic of $P(x,D)$ intersects $\Lambda_0$ a finite number of times.  Then, if $T \in I^{p,l}(\triangle, C_{P})$, $$T: I^r(\Lambda_0) \rightarrow I^{r+p, l}(\Lambda_0, \Lambda_1),$$ where $\Lambda_1 = C_P \circ \Lambda_0$ is the flowout from $\Lambda_0$ on $char(P)$. 
\end{prop}

\begin{prop} \label{charflow}
Suppose $\Lambda_1 \subset T^*(X) \backslash 0$ is a conic Lagrangian which is characteristic for $P$, meaning that $\Lambda_1 \subset char(P)$.  Then if $T$ is as above, then $$T: I^r(\Lambda_1) \rightarrow I^{r+p+ \frac{1}{2}}(\Lambda_1)$$ and therefore $$Q:I^r(\Lambda_1) \rightarrow I^{r+p-\frac{1}{2}}(\Lambda_1).$$

\end{prop}
\noindent Both of these propositions will be used to describe the scattering kernel when solving the direct problem.

Because of the geometry that arises when dealing with an approximation to the scattering kernel, it is necessary to understand the action of an FIO associated to $C_P$ on paired Lagrangian distributions with wavefront containing characteristic points.  The following theorem appears as Proposition 4.1 in \cite{GuiUhl81} in a more general form.  We give an alternate proof in order to emphasize the multiphase interpretation of paired Lagrangian distributions.

\begin{thm} \label{fiomapping}
Let $u \in I^{M_1, M_2}(\Lambda_1, \Lambda_2)$, with the Lagrangians intersecting in codimension $d$, and $F \in I^{p}(\Gamma)$, where $\Gamma$ is a homogeneous canonical relation such that $\Gamma \circ \Lambda_1 = \tilde{\Lambda}_1$ and $\Gamma \circ \Lambda_2 = \tilde{\Lambda}_2$ and the compositions are transversal.  If $\tilde{\Lambda}_1$ and  $\tilde{\Lambda}_2$ are cleanly intersecting in codimension $d$ as well, then 

\begin{equation}
Fu \in I^{M_1 + p, M_2}(\tilde{\Lambda}_1, \tilde{\Lambda}_2).
\end{equation}

\end{thm}

\noindent \underline{Proof}:  By a result in \cite{MelUhl79}, we can conjugate by elliptic Fourier integral operators associated with a canonical transformation and therefore assume that our domain is $\mathbb{R}^n$, $\Lambda_1=N^*(\{x'=0\})$, $\Lambda_2=N^*(\{x'=x''=0\})$, where $(x',x'') \in \mathbb{R}^{m+d}$.  The oscillatory representation of $u$ is now \begin{equation} \label{u}
\int e^{i[y'\cdot \xi' + y'' \cdot \xi'']} b(y; \xi', \xi'') d \xi' d \xi'',
\end{equation}
where $b \in S^{M_1,M-2}(\mathbb{R}^{n} \times(\mathbb{R}^{n-d} \times \mathbb{R}^d) \backslash 0)$.  For $F$, the Schwarz kernel takes the form $$\int e^{i\phi(x,y;\theta)} a(x; \theta) d \theta$$ with $\phi(x,y;\theta)$ parametrizing our (transformed) canonical relation $\Gamma$.  In order to show that $Fu(x)$ is another paired Lagrangian, we break our region of integration into subregions and show the new phase is a multiphase and the product of symbols satisfies a symbol-valued symbol estimate in each subregion.  The oscillatory representation of $Fu(x)$ is 
\begin{equation}\label{oscill}
\int e^{i[\phi(x,y; \theta) + y'\cdot \xi' + y'' \cdot \xi'']} a(x,y; \theta) b(y; \xi', \xi'') d \theta d \xi' d \xi'' dy.
\end{equation}
\noindent There are three regions to consider:  

\vskip.2cm \noindent \underline{Case 1}: If $\{|\xi'| \approx |\xi''|\}$, then $u \in I^{M_1 + M_2}(\Lambda_2)$.  The symbol $b$ in (\ref{oscill}) now satisfies a standard symbol estimate and we follow the proof of Hormander's result on the composition of Fourier integral operators whose canonical relations intersect transversally.  This requires us to define a new phase variable $\omega = ((|\theta|^2 +|\xi|^2)^{\frac{1}{2}}y, \theta, \xi))$. It follows that in the region $\{|\xi| \approx |\theta|\}$, our phase in (\ref{oscill}) parametrizes $\Gamma \circ \Lambda_2$ and the product of $a$ and $b$ satisfies a standard symbol estimate.  When $\{|\theta| \lesssim |\xi '| \approx |\xi ''|\}$, or vice versa, integration by parts in (\ref{oscill}) shows the order of the symbols can be decreased arbitrarily, making $Fu$ a $C^{\infty}$ function.  Hence, $Fu \in I^{M_1 + M_2 + p}(\tilde{\Lambda_2})$.

\vskip.2cm \noindent \underline{Case 2}: If $\{|\xi'| \lesssim |\xi''|\}$, then $u \in I^{M_1 + M_2}(\Lambda_2)$ once again and 
$Fu \in \in I^{M_1 + M_2 + p}(\tilde{\Lambda_2})$ by the above arguments.

\vskip.2cm \noindent \underline{Case 3}: Suppose $\{|\xi'| \gtrsim |\xi''|\}$.  In this region, we are microlocalized to the intersection $\Lambda_1 \cap \Lambda_2$.  First, if $\{y'' \neq 0\}$, then we can integrate by parts in (\ref{u}) in $\xi''$ arbitrarily many times to integrate out this variable and see that $u \in I^{M_1}(\Lambda_1)$.  Once again following H\"{o}rmander's argument and setting $\omega ' = ((|\theta|^2 +|\xi'|^2)^{\frac{1}{2}}y, \theta, \xi'))$, we get that $Fu \in I^{M_1 + p}(\tilde{\Lambda_1})$.

By our computations, it is natural to consider the new ``big" variable in the paired Lagrangian distribution as $\omega ' = ((|\theta|^2 +|\xi'|^2)^{\frac{1}{2}}y, \theta, \xi'))$ and the new ``small" variable as $\omega '' = \xi ''$.  Let us recall that on $\tilde{\Lambda_1}$, 
\begin{center}
$\begin{cases}
d_{y}\phi(x,y;\theta) + (\xi ',0) = 0 \\
d_{\theta}\phi(x,y;\theta) = 0 \\
y'=0,
\end{cases}$
\end{center}
and on $\Lambda_2$,
\begin{center}
$\begin{cases}
d_{y}\phi(x,y;\theta) + (\xi ',\xi '') = 0 \\
d_{\theta}\phi(x,y;\theta) = 0 \\
y'=y''=0.
\end{cases}$
\end{center}
Hence, in a conic neighborhood of the intersection, we are near the closed set
\begin{equation} \label{intersect}
\begin{cases} 
d_{y}\phi(x,y;\theta) + (\xi ',0) = 0 \\
d_{\theta}\phi(x,y;\theta) = 0 \\
y'=y''=0.
\end{cases}
\end{equation} 
As the phase function $\phi(x,y; \theta)$ is homogenous of degree 1 in the $\theta$ variable, we can assume $|d_{y}\phi| \approx |\theta|$.  It follows that within a conic neighborhood of (\ref{intersect}), $|\theta| \approx |\xi'|$.  The symbol in (\ref{oscill}), which we will refer to as $c(x; \omega ', \omega '')$, will now lie inside of $S^{M_1+p, M_2}$.  $\blacksquare$

\begin{rem} \label{intuition}
A corollary to the version of this theorem that appears in \cite{GuiUhl81} says if $\Gamma$ is a homogeneous canonical relation, satisfying the same transversality conditions with respect to the codimension $d$ intersecting pair of Lagrangians $(\Lambda_1, \Lambda_2)$, then $(\Gamma \circ \Lambda_1, \Gamma \circ \Lambda_2)$ is a codimension $d$ cleanly intersecting pair once again.  We could have used this theorem to do away with the clean intersection calculations in Section \ref{lscr}.  However, the computations themselves give us a deeper understanding of the present geometry, reappearing in calculations that show the original surfaces $S_1$ and $S_2$ are one-to-one with our flowouts and reflected Lagrangians, which is the crux of the inverse scattering problem.
\end{rem}

Let us now analyze the application of $\square^{-1}$ to a nested conormal distribution that will appear in a later section.  Using the oscillatory representation of $q\delta$ in $I^{M_1, M_2}(S_{1+}, S_{2+})$ discussed in Section \ref{mnc}, 
\begin{eqnarray} \label{bigintegral}
\nonumber \square^{-1}q\delta &=& \int e^{i[(x-y)\cdot \xi + (t-s)|\xi| + (t-s) \rho] + (s- y \cdot \omega) \tau + \vec{h}_1(y) \cdot \theta ' + \vec{h}_2(y) \cdot \theta '']} \\
&& \times a(x,y; \rho, \xi) b(y; \theta ', \theta '') \hskip.1cm  d \rho d \xi d \tau d \theta ' d \theta '' dy ds.
\end{eqnarray}
where $b$ is localized to the region $\{|\tau| \approx |\theta '| \gtrsim |\theta ''| \}$ in order capture the intersection $\Lambda_{1+} \cap \Lambda_{2+}$.

To microlocalize to $C_{\square} \cap \triangle$, we localize $a$ to the region $\{ |\xi| \gtrsim |\rho| \}$. If we do not focus on this part of phase space, then we are microlocalized to $\triangle$ away from $C_{\square}$ and are therefore applying a pseudodifferential operator; the wavefront set of $\square^{-1}q\delta(t-x \cdot \omega)$ does not move if $\square^{-1}$ acts like a pseudodifferential operator. Wavefront set calculus tells us 
\begin{equation}
WF(\square^{-1}q\delta(t-x \cdot \omega)) \subset \Lambda_{1+} \cup \Lambda_{2+} \cup \Lambda^a_{1+} \cup \Lambda^c_{2+}.
\end{equation}
The calculations in Section \ref{lscr} show that these Lagrangians all share a common submanifold of intersection, namely $\Sigma_1 \cup \Sigma_2$.  By interpreting the phase in (\ref{bigintegral}) as a generalized multiphase, integration by parts and application of the standard theory of Fourier integral operators in (\ref{bigintegral}) shows that these singularities do appear.  A theory of intersecting quadruples of Lagrangians in a certain configuration does appear in \cite{MelUhl79} but does not apply in this setting because certain conditions involving the arrangement of the Lagrangians are not met. 

If we further localize our symbols to $\{|\rho| \approx |\tau| \}$, then letting our ``big" variable be $((|\xi|^2 + |\tau|^2 + |\theta '|^2 )^{\frac{1}{2}}(y,s),\xi, \tau, \theta ')$, our ``medium" variable be $\rho$, and our ``small" variable be $\theta ''$, we get a distribution whose wavefront set is now inside of $\Lambda^a_{1+} \cup \Lambda_{1+} \cup \Lambda_{2+}$.  We will elaborate further on this triple in Section \ref{bose}.

The final microlocalization will be to further localize our symbols to the region $\{|\xi| \approx |\tau| \}$ in phase space and to $\{t \neq s\}$ in space.  Integrate by parts arbitrarily many times in $\rho$ to obtain a new symbol $c(x,y; \xi, \tau, \theta ', \theta '') = b(y; \theta ', \theta '') \times \int e^{i[(t-s) \rho]} a(x,y; \xi, \rho) \hskip.1cm d \rho$ satisfying a symbol-valued symbol estimate 
\begin{equation}
\partial_{x,y}^{\gamma} \partial_{\theta ''}^{\beta} \partial_{\xi, \tau, \theta '}^{\alpha}c(x,y; \xi, \tau, \theta ', \theta '') \lesssim \langle \xi, \tau, \theta ', \theta '' \rangle ^{M_1-\frac{1}{2} - |\alpha|} \langle \theta '' \rangle ^{M_2 - |\beta|}.
\end{equation} 
We now let our ``big" variable be $((|\xi|^2 + |\tau|^2 + |\theta '|^2 )^{\frac{1}{2}}(y,s),\xi, \tau, \theta ')$ and ``small" variable be $\theta ''$.  This stratification puts the distribution into the regime of a paired Lagrangian distribution associated to the codimension $d_2$ cleanly intersecting pair of $(\Lambda^a_{1+}, \Lambda^c_{2+})$.  Section \ref{sip} shows that this portion of $\square^{-1}q \delta$ contains all the information needed to solve the inverse problem.

\subsection{Sobolev estimates}

\noindent We equip $T^*(\mathbb{R}^n) \backslash 0$ with the canonical symplectic two-form $\omega = \Sigma d\xi_i \wedge d x_i$.  Let $\Sigma \subset T^*(\mathbb{R}^n)$ is a smooth, codimension $k$ conic submanifold with $1 \leq k < n$ that is also involutive, meaning that the ideal of smooth functions that vanish on $\Sigma$ is closed under the Poisson bracket.  Theorems in symplectic geometry from Section 3 of \cite{Duis11} show that the flowout of $\Sigma $, $\Lambda_{\Sigma} \in T^*(\mathbb{R}^n) \backslash 0 \times T^*(\mathbb{R}^n) \backslash 0$, is a canonical relation given by 
\begin{equation}
C_{\Sigma} = \{\big((x;\xi),(y;\eta)\big)\in \Sigma \times \Sigma : (y; \eta) \in \Xi_{(x; \xi)}\},
\end{equation}
where $\Xi_{(x; \xi)}$ is the bicharacteristic leaf of $\Sigma$ containing $(x;\xi)$.  It is straightforward to see $C_{\Sigma}$ and $\triangle$ intersect in codimension $k$; refer to Section 4 of \cite{Duis11}.  The following is a theorem from \cite{GU90}.

\begin{thm} \label{boxsobo}
Let $A \in I^{p,l}(\mathbb{R}^n \times \mathbb{R}^n; \triangle, C_{\Sigma})$.  Then $$A: H^s_{comp}(\mathbb{R}^n) \rightarrow H^{s+s_0}(\mathbb{R}^n)$$ continuously for all $s \in \mathbb{R}$ if $$\sup\{p + \frac{k}{2}, p + l\} \leq -s_0.$$ 
\end{thm}

\noindent Consider the flowout relation for $\square^{-1}$ in (\ref{flowout}).  The parametrization shows that the cotangent variables to $(x,t)$ sit in a hypersurface in the fibers of our phase space, namely the light cone.  Because this submanifold has codimension 1, it follows by a rank calculation that $k=1$ in the above theorem for $A=\square^{-1}$.  Hence for $-s_0 \geq \text{max}(-1,-2) = -1$, $\square^{-1}$ is smoothing of order at least 1.  In fact, away from the characteristic variety, $\square^{-1}$ acts as a pseudodifferential operator of order $-2$ and increases Sobolev regularity by $2$.  

\begin{rem}
It turns out that multiplication by standard conormal distributions fall under the hypothesis of this theorem.
\end{rem}

\section{Born series} \label{bose}

\noindent We analyze the formal series 
\begin{equation} \label{born approx}
\sum_{i}^{\infty}{(-1)^i(\square^{-1}M_q)(\delta(t-x\cdot \omega))}, 
\end{equation} 
which in the physics literature is called the \underline{Born series}, in order to solve the direct problem.  Set $u_i:=(-1)^i(\square^{-1} M_q)(\delta(t-x \cdot \omega))$ where $u_0=\delta(t-x\cdot \omega)$ so that the series on the right of (\ref{born approx}) is formally telescoping when $\square + q$ is applied.  The remaining parts of this paper will show the first two terms,
\begin{equation*}
u_0 + u_1 = \delta(t - x\cdot \omega) - (\square^{-1}M_q)(\delta(t-x\cdot \omega)),
\end{equation*}
commonly referred to as the (first order) \underline{Born approximation}, will be enough to solve the inverse problem.  Analysis of the higher terms is substantially more intricate, involving a more complete theory of distributions associated to higher order systems of Lagrangians.  At this point, we do not know if an in-depth analysis is even possible without strict assumptions on the geometry of $S_1$ and $S_2$. 

\subsection{Computing the series}

\noindent For $q \in I^{M_1,M_2}(S_1,S_2)$, we apply a microlocal partition of unity $\sum_{i=1}^3 $ ${\chi_i}(x,t, \omega)$ that gives 
\begin{eqnarray}
\nonumber \chi_1 q &\in& I^{M_1}(S_1), \\
\nonumber \chi_2 q &\in& I^{M_1 + M_2}(S_2), \text{ and } \\
\chi_3 q &\in& I^{M_1,M_2}(S_1, S_2),
\end{eqnarray}
where the last term is microlocally supported near the intersection $\Lambda_1 \cap \Lambda_2$.  This partition is introduced to facilitate the analysis of the singularities of $u_1$.

By Lemmas \ref{singlemult} and \ref{doublemult} and remembering that supp $\delta \subset \{t-x \cdot \omega=0\}$, 
\begin{eqnarray} \label{u.5}
\nonumber \chi_1 q \cdot \delta(t - x \cdot \omega) &\in& I^{0,M_1}(S_+, S_{1+}),\\
\nonumber \chi_2 q \cdot \delta(t - x \cdot \omega) &\in& I^{0, M_1 + M_2}(S_+, S_{2+}), \text{ and } \\ 
\chi_3 q \cdot \delta(t - x \cdot \omega) &\in& I^{0,M_1,M_2}(S_+, S_{1+}, S_{2+}) + I^{M_1,M_2}(S_{1+}, S_{2+}).
\end{eqnarray}

\noindent Keeping in mind the ``good" part of $\Sigma^1 \cup \Sigma^2$ does not meet $\Lambda_{+} \cap \Lambda_{1+}$ and $\Lambda{+} \cap \Lambda_{2+}$ , or rather  $\Lambda_{+} \cap \Lambda_{1+} \subset \Sigma^1 \backslash \Sigma_1$ and $\Lambda_{+} \cap \Lambda_{2+} \subset \Sigma^2 \backslash \Sigma_2$, we can microlocalize away from a conic neighborhood $\mathcal{O}$ of $(\Sigma^1 \backslash \Sigma_1) \cup (\Sigma^2 \backslash \Sigma_2)$, changing (\ref{u.5}) to
\begin{eqnarray} \label{u.55}
\nonumber && \chi_1 q \cdot \delta(t - x \cdot \omega) \in I^{\frac{1}{2} - \frac{n}{2}}(\Lambda_{+} \backslash (\mathcal{O} \cap \Lambda_+)) + I^{M_1 + \frac{d_1+1}{2} - \frac{n}{2}}(\Lambda_{1+} \backslash (\mathcal{O} \cap \Lambda_{1+})) \\
\nonumber &&\indent + H^{s_0}_{\mathcal{O}}, \\
\nonumber && \chi_2 q \cdot \delta(t - x \cdot \omega) \in I^{\frac{1}{2} - \frac{n}{2}}(\Lambda_+ \backslash (\mathcal{O} \cap \Lambda_+))  \\
\nonumber && \indent + I^{M_1 + M_2 + \frac{d_1+d_2+1}{2} - \frac{n}{2}}(\Lambda_{2+} \backslash (\mathcal{O} \cap \Lambda_{2+}))+ H^{s_0}_{\mathcal{O}}, \\
\nonumber && \chi_3 q \cdot \delta(t - x \cdot \omega) \in I^{\frac{1}{2} - \frac{n}{2}}(\Lambda_{+} \backslash (\mathcal{O} \cap \Lambda_+)) \\
\nonumber && \indent + I^{M_1+M_2+\frac{d_1+d_2+1}{2} - \frac{n}{2},-M_2-\frac{d_2}{2}}(\Lambda_{1+} \backslash (\mathcal{O} \cap \Lambda_{1+}), \Lambda_{2+} \backslash (\mathcal{O} \cap \Lambda_{2+}))\\
&& \indent + H^{s_0}_{\mathcal{O}}, 
\end{eqnarray}
where we have switched from conormal notation to Lagrangian notation and $H^{s_0}_{\mathcal{O}}$ are elements $w \in H^{s_0}(\mathbb{R}^n \times \mathbb{R} \times S^{n-1})$ with $WF(w) \subset \mathcal{O}$.  Note the sums of spaces are not direct sums as there is some overlap in the singularities from one term to the next.

We are now ready to describe the second term, $u_1$, of the Born approximation.  Let $ L := C_{\square} \circ \mathcal{O}$, where $C_{\square}$ is the flowout relation for $\square^{-1}$.  The set $L$ will be invariant under the Hamiltonian flow of $H_{\square}$.  Applying $\square^{-1}$ to both sides of (\ref{u.55}), using Propositions \ref{singleflow} and \ref{charflow} and Theorem \ref{fiomapping} gives 
\begin{eqnarray} \label{u_1}
\nonumber u_1 = \square^{-1}q\delta &\in& I^{\frac{1-n}{2} - 1}(\Lambda_+ \backslash L_+) + I^{M_1 + \frac{d_1 + 1}{2} - \frac{n}{2} - \frac{3}{2}, -\frac{1}{2}}(\Lambda_{1+} \backslash L_{1+}, \Lambda_{1+}^a \backslash L^a_{1+})\\
\nonumber &+&  I^{M_1 + M_2 + \frac{d_1 + d_2 + 1}{2} - \frac{n}{2} - \frac{3}{2}, -\frac{1}{2}}(\Lambda_{2+} \backslash L_{2+}, \Lambda_{2+}^c \backslash L_{2+}^c)\\
 &+&  I(\Lambda_{1+}^a \backslash L_{1+}^a, \Lambda_{2+}^c \backslash L_{2+}^2, \Lambda_{1+} \backslash L_{1+}, \Lambda_{2+} \backslash L_{2+}) + H_{L}^{s_0+1}
\end{eqnarray}
where the last term $I(\Lambda_{1+}^a \backslash L_{1+}^a, \Lambda_{2+}^c \backslash L_{2+}^2, \Lambda_{1+} \backslash L_{1+}, \Lambda_{2+} \backslash L_{2+})$ in (\ref{u_1}) is a new class of distributions such that $WF(u) \subset \Lambda_{1+}^a \cup \Lambda_{2+}^c \cup \Lambda_{1+} \cup \Lambda_{2+}$, $H_{L}^{s_0+1}$ are Sobolev elements of order $s_0+1$ with wavefront set inside of $L$, and the varied labeling of $L$ denotes intersection of $L$ with the juxtaposed Lagrangian.  Note that we have used Theorem \ref{boxsobo} on the Sobolev space term.  The quadruple portion of $u_1$ was analyzed earlier in Section \ref{pwo}.

In Section \ref{pwo}, we prove, microlocally away from $\Lambda_{1+}^a \cap \Lambda_{2+}^c \cap \Lambda_{1+} \cap \Lambda_{2+}$, there exists a nontrivial singularity on $\Lambda_{1+}^a \cap \Lambda_{1+} \cap \Lambda_{2+}$.  We label this class similarly as $I(\Lambda_{1+}^a, \Lambda_{1+}, \Lambda_{2+})$, the space of distributions whose wavefront set is contained inside the union of the three listed Lagrangians.  Notice we do not attach orders to this newly defined class in spite of it being similarly defined to those of the nested triple conormal configuration.  This is because for codimensions of $S_1$ greater than 1, there exists a conic singularity on $S_{1+}^a$ along $S_{1+}$.  Hence, because of the degeneracy of the underlying submanifolds, this new triple does not fall into our previously defined class.  However, for $d_1=1$, the distribution is a nested triple conormal of order $(-\frac{1}{2}, M_1, M_2)$.  The space
\begin{equation}
I^{M_1 + M_2 + \frac{d_1 + d_2 + 1}{2} - \frac{n}{2} - \frac{1}{2}, -M_2 - \frac{d_2}{2}}(\Lambda_{1+}^a, \Lambda_{2+}^c)
\end{equation}
also appear as a portion in the description of the quadruple term from (\ref{u_1}).  After all these observations, it follows that 
\begin{eqnarray}
\nonumber u_1 &\in& I^{\frac{1-n}{2} - 1}(\Lambda_+ \backslash L_+) + I^{M_1 + \frac{d_1 + 1}{2} - \frac{n}{2} - \frac{3}{2}, -\frac{1}{2}}(\Lambda_{1+} \backslash L_{1+}, \Lambda_{1+}^a \backslash L_{1+}^a)\\
\nonumber &+&  I^{M_1 + M_2 + \frac{d_1 + d_2 + 1}{2} - \frac{n}{2} - \frac{3}{2}, -\frac{1}{2}}(\Lambda_{2+}\backslash L_{2+}, \Lambda_{2+}^c \backslash L_{2+}^c) \\
\nonumber &+& I_{\text{loc}}(\Lambda_{1+}^a\backslash L, \Lambda_{2+}^c\backslash L, \Lambda_{1+}\backslash L_{1+}, \Lambda_{2+}\backslash L_{2+})\\
\nonumber &+& I_{\text{loc}}(\Lambda_{1+}^a\backslash L_{1+}^a, \Lambda_{1+}\backslash L_{1+}, \Lambda_{2+}\backslash L_{2+})\\
&+& I^{M_1 + M_2 + \frac{d_1 + d_2 + 1}{2} - \frac{n}{2} - \frac{1}{2}, -M_2 - \frac{d_2}{2}}(\Lambda_{1+}^a\backslash L_{1+}^a, \Lambda_{2+}^c\backslash L_{2+}^c) + H_{L}^{s_0+1}
\end{eqnarray}
where the subscript ``loc" denotes localization near the respective triple and quadruple intersections.  As the $t$ variable is bounded on $\Lambda_{1+}$ and $\Lambda_{2+}$ because the support of $u_1$ is inside $\{t-x \cdot \omega=0\}$, for $t >> 0$,
\begin{eqnarray} 
\nonumber u_1 &\in& I^{\frac{1-n}{2} - 1}(\Lambda_+ \backslash L_+) + I^{M_1 + \frac{d_1 + 1}{2} - \frac{n}{2} - \frac{3}{2}}(\Lambda_{1+}^a \backslash L_{1+}^a) \\
\nonumber &+& I^{M_1 + M_2 + \frac{d_1 + d_2 + 1}{2} - \frac{n}{2} - \frac{3}{2}}(\Lambda_{2+}^c \backslash L_{2+}^c) \\
&+& I^{M_1 + M_2 + \frac{d_1 + d_2 + 1}{2} - \frac{n}{2} - \frac{1}{2}, -M_2 - \frac{d_2}{2}}(\Lambda_{1+}^a \backslash L_{1+}^a, \Lambda_{2+}^c \backslash L_{2+}^c) + H_{L}^{s_0+1}.
\end{eqnarray}

\subsection{Subsequent terms and singularities} 

\noindent We show in this section that the remaining part of the Born series, i.e. $\bar{u}=u_2 + u_3 + ...$, has $WF$ which grows larger with each iteration of $\square^{-1}M_q$.  As a consequence, we will use the mapping properties of $\square^{-1}$ and $M_q$ from previous sections to place the $\bar{u}$ into a suitable Sobolev of higher regularity than the previous two terms.  This final step will be taken in Section \ref{sip}.

It was shown in Section \ref{lscr} that $\Lambda_{1+}^a$ and $\Lambda_{2+}^c$ are conormal bundles of the hypersurfaces $S_{1+}^a$ and $S_{2+}^c$, respectively.  A well-known fact about the multiplication of distributions $u,v \in \mathcal{D}'(X)$, where $X$ is an open set in $\mathbb{R}^{n}$, states if $WF(u) \cap WF'(v) = \emptyset$ then 
\begin{equation*}
WF(uv) \subset WF(u) + WF(v).
\end{equation*}
where $WF(u) + WF(v) = \{(x, \xi_u + \xi_v) : (x, \xi_u) \in WF(u), (x, \xi_v) \in WF(v)\} \subset T^*(\mathbb{R}^n \times \mathbb{R} \times S^{n-1}) \backslash 0$; see \cite{Ho71}.

If we microlocalize $u \in I^{M_1 + M_2 + \frac{d_1 + d_2 + 1}{2} - \frac{n}{2} - \frac{3}{2}, -\frac{1}{2}}(\Lambda_{2+}, \Lambda_{2+}^c)$ away from $\Lambda_{2+}$, we get $u \in I^{M_1 + M_2 + \frac{d_1 + d_2 + 1}{2} - \frac{n}{2} - \frac{3}{2}}(\Lambda_{2+}^c)$.  In fact, $u$ is conormal to $S_{2+}^c$ away from $\Lambda_{2+} \cap \Lambda_{2+}^c$.  It is the interaction of this part of $u_1$ with  $\chi q \in I^{M_1}(S_1)$ that creates the new $WF$.

Using Lemma \ref{singlemult}, all we must show is that there exist submanifolds $S_1$ and $S_2$ such that $S_1 \pitchfork S_{2+}^c$ and this intersection is different from all other submanifolds of the same codimension.

\begin{prop} 
Let $S_1 = \{x_n=0\}$ and $S_2=\{x_n=x_{n-1}=0\}$ be subsets of $\mathbb{R}^{n+1} \times \mathbb{S}^{n-1}$.  Then $S_1 \pitchfork S_{2+}^c$ is a new submanifold of codimension $2$ in $\mathbb{R}^{n+1} \times \mathbb{S}^{n-1}$ that is different from those that $u_1$ is conormal.
\end{prop}

\noindent \underline{Proof:} Considering the general parametrization of $\Lambda_{2+}^c$ given in Section 5, we apply the Implicit Function Theorem to solve for the $y_{n-2}$ coordinate.  More explicitly,
\begin{eqnarray*}
S_{2+}^c &=& \{((y'',0,0) - r(\theta_1 \vec{e}_{n-1} + \theta_2 \vec{e}_{n} - \frac{\theta_1^2 + \theta_2^2}{2(\theta_1 \omega_{n-1} + \theta_2 \omega_{n})}\omega), y^{''} \cdot \omega^{''} \\ 
&& + r\frac{\theta_1^2 + \theta_2^2}{2(\theta_1 \omega_{n-1} + \theta_2 \omega_{n})}, \omega) :  y_{n-2}=f(y_1,...,y_{n-3}), (\theta_1, \theta_2) \in \mathbb{R}^2 \backslash 0, \\
&& \omega \in \mathbb{S}^{n-1}, r \in \mathbb{R}\}
\end{eqnarray*}
\noindent From the parametrization of $S_1$, which is 
\begin{equation*}
\{ ((y'',y_{n-1},0),s, \omega) : (y'', y_{n-1}) \in \mathbb{R}^{n-1}, s \in \mathbb{R}, \omega \in  \mathbb{S}^{n-1} \},
\end{equation*} 
it follows $S_3 = S_1 \cap S_{2+}^c$ has the same parametrization as $S_{2+}^c$ but with $r = \frac{\theta_1^2 + \theta_2^2}{2(\theta_1 \omega_{n-1} + \theta_2 \omega_{n})} \cdot \frac{1}{\theta_2} \cdot \omega_n$.  For brevity in the following argument, set $\tau = \frac{\theta_1^2 + \theta_2^2}{2(\theta_1 \omega_{n-1} + \theta_2 \omega_{n})}$

Computing tangent vectors, we see $\{ \{ \frac{\partial}{\partial x_i}\}_{i=1}^{n-1}, \frac{\partial}{\partial s}, \{ \frac{\partial}{\partial \omega_i}\}_{i=1}^{n} \}$ on $S_1$ spans a $2n-1$ dimensional set in $T_x(\mathbb{R}^{n+1} \times \mathbb{S}^{n-1})$ for $x \in S_3$.  Take $\frac{\partial}{\partial \theta_2}$ on $S_{2+}^c$ and observe its $n$th spatial coordinate, which is $\frac{\tau}{\theta_2}(1 - \frac{\partial}{\partial \theta_2} \tau \omega_n).$  This quantity is not identically equal to $0$ because 0-sections are deleted from the total space and the ability to work away from $(\frac{\partial}{\partial \theta_2}\tau) \omega_n = 1$, a lower-dimensional set that is independent of $r$.  Hence, at $x \in S_3$, $T_x(S_1) \oplus T_x(S_{2+}^c) = T_x(\mathbb{R}^{n+1} \times \mathbb{S}^{n-1})$ and the intersection is transverse.

We now show $S_3$ is not one of the previously encountered submanifolds.  It is clear $S_2$ has its $x_{n-1}$ coordinate identically equal to 0.  However, $\frac{\tau}{\theta_2}(\theta_1 - \tau \omega_{n-1}) \neq 0$, almost everywhere.  This shows $S_2 \cap S_3 \subsetneq S_3$.  In $S_{1+}$, the relation $t - x \cdot \omega = 0$ is satisfied.  When this relation is applied to the parametrization of $S_3$, it follows that $-r(\theta_1 \omega_{n-1} + \theta_2 \omega_n)=0$.  This is another equation that holds on a lower dimensional set, proving $S_{1+} \cap S_3 \subsetneq S_3$. $\blacksquare$

It follows from Lemma \ref{singlemult} that $q(x) \cdot \square^{-1}(q(x) \cdot \delta(t - x \cdot \omega)$ has different wave front set than does $u_0 + u_1$.  This gives the evidence, after application of $\square^{-1}$ to the above expression, that the new flowouts generated are different from $\Lambda_{1+}^a$ and $\Lambda_{2+}^c$, as $N^*(S_3)$ has characteristic points because $S_3 \subset S^c_{2+}$.  Hence, the wave front set of the Born series does not stabilize.  We will show in the next section that the approximate solution $u_0 + u_1$ will be enough to solve the inverse problem by simply including the remaining terms of the Born series in a Sobolev space.

\section{Solution to the inverse problem} \label{sip}

\noindent We now solve the inverse scattering problem by combining all of the results from previous Sections.  First, we define the scattering kernel, which describes the far-field pattern of a solution $u$ to the direct problem, using the theory of Lax and Phillips \cite{LP89}.  We then show that the restriction of the scattering kernel, known as the backscattering, also determines the singularities of $q$.  In fact, we prove a corollary that shows for well-behaved subsets of scattering data $\mathbb{D}$, the corresponding restriction of the scattering kernel continues to determine the singularities of $q$.

\subsection{Lax-Phillips theory and the scattering kernel}

\noindent The Lax-Phillips Radon transform \cite{LP89} is the map $R_{LP}: \mathcal{D}'(\mathbb{C}^2) \rightarrow \mathcal{D}'(\mathbb{C})$ and is defined by
\begin{equation}
R_{LP}(v_0,v_1)=C_n|D_s^{\frac{n-1}{2}}|(D_sR(v_0)-R(v_1))
\end{equation}
where $C_n$ is a constant depending on $n$, $|D_s^{\frac{n-1}{2}}|$ is a pseudodifferential operator acting on the $s$ variable of order $\frac{n-1}{2}$ which overlaps with a differential operator for $n$ odd, and $R$ is the Radon transform defined in Section \ref{lscr}.  To define the scattering kernel, we set $w=|D_t^{\frac{n-3}{2}}|(u(x,t,\omega) - \delta(t-x \cdot \omega))$ and get the relation
\begin{equation}
\alpha_q(t-s, \phi, \omega) = \delta(t-s) \otimes \delta(\phi - \omega) + R_{LP}(w, D_tw), \text{ for } t>>0.
\end{equation}
Focusing our attention on the scattered part of $u$ which is expressed as $u - \delta$, it follows that 
\begin{equation}\label{semiscatt}
(\alpha_q-\delta \otimes \delta)(t-s, \phi, \omega)=F(u-\delta)
\end{equation}
where $F \in I^{\frac{n-1}{2}}(C_R)$ and $C_R$ is the canonical relation associated to the Radon transform; see Section \ref{lscr}.  Since $F$ is elliptic on the first components of $C_R \circ \Lambda_{1+}^a$ and $C_R \circ \Lambda_{2+}^c$, these are the only sets that stay after the application of $C_{\rho}$, also defined in Section \ref{lscr}.  The reason for introducing $\rho^*$ is to deal with the translation invariance of (\ref{semiscatt}) for $t>>0$.  Summing all this up, we get 

\begin{defn}
The \underline{scattering kernel} associated with $q$ is $\alpha_q$ defined by
\begin{equation}
(\alpha_q-\delta \otimes \delta)(s, \phi, \omega)=\rho^*F(u-\delta),
\end{equation}
where $u$ is a solution to the continuation problem.
\end{defn}

\noindent Unlike in the case for scattering theory for the wave equation when $q \in C^{\infty}_0(\mathbb{R}^n)$ where there is a wave group \cite{MelUhlnotes}, we use the Born series to provide information about the solution $u$. It is natural to consider the \underline{approximate} \underline{scattering kernel} $\rho^*F(u_0 + u_1 - \delta)=\rho^*F(u_1)$ and investigate the error of the difference between this and the true scattering kernel.  The next section will show this error term is of lower order, giving us the ability to work with the approximate scattering kernel instead of the exact one.

\subsection{Comparison of scattering kernels}

\noindent Our principal objective is to make a comparison between $\alpha_q$ and $\rho^*F(u_0 + u_1 - \delta)$, the exact scattering kernel and the approximate scattering kernel.  All of the analysis from our previous sections and further computations below will show that this difference of scattering kernels is negligible in an appropriate sense.

For now, set $\bar{u} := u - (u_0 + u_1)$.  For $t>>0$, 
\begin{eqnarray} \label{setup}
\nonumber (\square + q)\bar{u} &=& 0 - \square u_0 - q \cdot u_0 - \square u_1 - q \cdot u_1\\
&=& - q \cdot u_1
\end{eqnarray}
We recall that for the same range of $t$,
\begin{eqnarray}\label{spaces}
\nonumber u_1 &\in& I^{\frac{1-n}{2} - 1}(\Lambda_+ \backslash L_+) + I^{M_1 + \frac{d_1 + 1}{2} - \frac{n}{2} - \frac{3}{2}}(\Lambda_{1+}^a \backslash L_{1+}^a) \\
\nonumber &+& I^{M_1 + M_2 + \frac{d_1 + d_2 + 1}{2} - \frac{n}{2} - \frac{3}{2}}(\Lambda_{2+}^c \backslash L_{2+}^c) \\
&+& I^{M_1 + M_2 + \frac{d_1 + d_2 + 1}{2} - \frac{n}{2} - \frac{1}{2}, -M_2 - \frac{d_2}{2}}(\Lambda_{1+}^a \backslash L_{1+}^a, \Lambda_{2+}^c \backslash L_{2+}^c) + H_{L}^{s_0+1}.
\end{eqnarray}
Moreover, 
\begin{eqnarray}\label{diff}
\nonumber \rho^*F(u - \delta) - \rho^*F(u_0 + u_1 - \delta) &=& \rho^*F(u - (u_0 + u_1))\\
&=& \rho^*F(\bar{u}), 
\end{eqnarray}
showing why placing $\rho^*F(\bar{u})$ into a space of lower order will allow us to extract information from the true scattering kernel by reading the approximate scattering kernel.  It is necessary to prove this error term does not interfere by computing it's Sobolev regularity.

We define  
\begin{equation} \label{Neumser}
(I + \square^{-1}M_q)^{-1} := \sum_{j=0}^{\infty}(-1)^j (\square^{-1}M_q)^j,
\end{equation}
as this makes the Born series consistent with a Neumann series; also, we only need to know the Born series asymptotically.  In addition, $\square^{-1}M_q$ must raise regularity on Sobolev spaces in order to have the subsequent terms of the Born series be smoother.  This happens when the orders of $q$ satisfy $M_1+\frac{M_2}{2} < -d_1 - d_2 + 1$ for the specific ranges of $M_1$ and $M_2$ discussed in Theorem \ref{sobomult}.  Recall that this allows $q$ to blowup on parts of $S_1$ or $S_2$ depending on the specific orders.  Applying $\square^{-1}$ to both sides of (\ref{setup}) gives the new relation 
\begin{equation}\label{setup2}
(I+ \square^{-1}M_q)\bar{u}= u_2 := \square^{-1}M_q(u_1).
\end{equation}
Because of the iterated regularity characterizations of the spaces that appear in (\ref{spaces}), we can place $u_1$ into a Sobolev space of some order $\tilde{s}$.  Another application of $\sum_{j=0}^{N-1}(-1)^j (\square^{-1}M_q)^j$ to both sides of (\ref{setup2}) leads to 
\begin{equation} \label{blah}
(-1)^{N-1}(\square^{-1}M_q)^N\bar{u} = \sum_{j=0}^{N-1}(-1)^j (\square^{-1}M_q)^j( u_2),
\end{equation}  
with $\square^{-1}M_q$ raising Sobolev regularity by some amount in each term.  It follows that the right hand side of (\ref{blah}) stays in the Sobolev space $H^{\tilde{s}}$ by convergence of the Neumann series.  Hence $\bar{u} \in H^{\tilde{s}}$.  Set $\hat{L} := C_{\rho} \circ C_{R} \circ L$.  Keeping in mind the mapping properties of $\rho^*$ and $F$ on Sobolev spaces \cite{Foll95, Ho71} and spaces of Lagrangian distributions from the transverse intersection calculus discussed in Section \ref{pre}, (\ref{spaces}) implies 

\begin{thm}\label{comparison}
The difference of the scattering kernels, $\alpha(s,\phi, \omega) - \rho^*F(u_0 + u_1 - \delta)$, is in $H^{\tilde{s}}$ and
\begin{eqnarray}
\alpha(s, \phi, \omega) &\in& I^{\frac{1-n}{2} - 1 + \frac{3-2n}{4}}(\hat{\Lambda}_+ \backslash \hat{L}_+) \\
\nonumber &+& I^{M_1 + \frac{d_1 + 1}{2} - \frac{n}{2} - \frac{3}{2} + \frac{3-2n}{4}}(\hat{\Lambda}_{1+}^a \backslash \hat{L}_{1+}^a) \\
\nonumber &+& I^{M_1 + M_2 + \frac{d_1 + d_2 + 1}{2} - \frac{n}{2} - \frac{3}{2} + \frac{3-2n}{4}}(\hat{\Lambda}_{2+}^c \backslash \hat{L}_{2+}^c) \\
\nonumber &+& I^{M_1 + M_2 + \frac{d_1 + d_2 + 1}{2} - \frac{n}{2} - \frac{1}{2} + \frac{3-2n}{4}, -M_2 - \frac{d_2}{2}}(\hat{\Lambda}_{1+}^a \backslash \hat{L}_{1+}^a, \hat{\Lambda}_{2+}^c \backslash \hat{L}_{2+}^c)\\
\nonumber &+& H_{\hat{L}}^{s_0 + \frac{3-2n}{4}} + H^{\tilde{s}}.
\end{eqnarray}
Hence, the principal symbol of $\alpha$ is same as that of $\rho^*F(u_0 + u_1 - \delta)=\rho^*F(u_1)$.
\end{thm}

\subsection{Determination of $S_1$, $S_2$, and $\mu(q)$}

\noindent A majority of this final section is essentially an adaptation of the conclusion of \cite{GU93} to the present context, with some minor modifications.  In order to even retrieve the singularities of $q$ from the restriction of $\alpha(s, \phi, \omega)$ to various sets of scattering data, i.e. submanifolds of $\mathbb{R} \times S^{n-1} \times S^{n-1}$, we must be able to do this without any restriction, i.e. just from using the full scattering kernel.  By Theorem \ref{comparison} of this section, it is enough to show $S_1$ and $S_2$ are determined by the reflected Lagrangians.  We will compare two reflected Lagrangians:
\begin{eqnarray}
\nonumber \hat{\Lambda}_{1+}^a &=& \{(-y \cdot (\frac{(\nu - \sigma \omega)}{\sigma} + \omega),- \frac{(\nu - \sigma \omega)}{\sigma}, \omega; \sigma, -\sigma i^*_{\frac{(\nu - \sigma \omega)}{\sigma}}(y), \sigma i^*_{\omega}(y)): \\
\label{reflect1} && (y, \nu) \in N^*(S_1), \omega \in S^{n-1}, r \in \mathbb{R}, \sigma \in \mathbb{R} \backslash 0\}
\end{eqnarray}
and
\begin{eqnarray}
\nonumber \bar{\hat{\Lambda}}_{1+}^a &=& \{(-\bar{y} \cdot (\frac{(\bar{\nu} - \bar{\sigma} \bar{\omega})}{\bar{\sigma}} + \omega),- \frac{(\bar{\nu} - \bar{\sigma} \bar{\omega})}{\bar{\sigma}}, \bar{\omega}; \bar{\sigma}, -\bar{\sigma} i^*_{\frac{(\bar{\nu} - \bar{\sigma} \bar{\omega})}{\bar{\sigma}}}(\bar{y}), \bar{\sigma} i^*_{\omega}(\bar{y})): \\
\label{reflect2} && (\bar{y}, \bar{\nu}) \in N^*(S_1), \bar{\omega} \in S^{n-1}, \bar{r} \in \mathbb{R}, \bar{\sigma} = \frac{|\bar{\nu}|^2|}{2 \bar{\nu} \cdot \bar{\omega}} \text{ with } \bar{\nu} \cdot \bar{\omega} \neq 0\}.
\end{eqnarray}
First noticing that $\omega=\bar{\omega}$ and $\sigma=\bar{\sigma}$, we set the $\phi$ coordinates equal to get $\nu=\bar{\nu}$.  The comparison of the $\Omega$ coordinates tells us $\bar{y}=y+c(y,\omega)\omega$. Now, set the $s$ coordinates equal to get the relation 
\begin{eqnarray}
\nonumber && -y \cdot (\frac{\nu}{\sigma}) = -\bar{y} \cdot (\frac{\nu}{\sigma}) \\
\nonumber && \Leftrightarrow - y \cdot \nu = -\bar{y} \cdot \nu \Leftrightarrow y \cdot \nu = (y+c(y,\omega)\omega) \cdot \nu \\
&& \Rightarrow c(y,\omega)\omega \cdot \nu = 0
\end{eqnarray}
Remembering $\omega \cdot \nu \neq 0$ because of the tangential rays condition, we get $c_1(y, \omega)=0$.  Therefore $y = \bar{y}$.  Hence $\hat{\Lambda}_{1+}^a$ determines $S_1$.  The same holds true for $\hat{\Lambda}_{2+}^c$ and $S_2$.

Let $\mathbb{B} = \{(s,\phi, \omega) \in \mathbb{R} \times S^{n-1} \times S^{n-1} : \phi=-\omega\}$ be the backscattering surface.  The map 
\begin{equation}
j_{\mathbb{B}}: \mathbb{R} \times S^{n-1} \rightarrow \mathbb{B}
\end{equation}
defined by $j_{\mathbb{B}}(s,\omega)=(s, \omega, -\omega)$ induces the pullback
\begin{equation}
j_{\mathbb{B}}^*: \mathcal{D}_{\mathbb{B}}'( \mathbb{R} \times S^{n-1} \times S^{n-1}) \rightarrow \mathcal{D}'(\mathbb{R} \times S^{n-1})
\end{equation}
with the domain being distributions whose wavefront set is disjoint from the normals of $j_{\mathbb{B}}$.  This pullback is another Fourier integral operator, $j_{\mathbb{B}}^* \in I^{\frac{n-1}{4}}(C_{\mathbb{B}})$ with 
\begin{eqnarray}
\nonumber C_{\mathbb{B}} &=& \{(s, \omega, \tau, \Omega; s', \phi, \omega ', \tau', \Phi, \Omega '):\\
\nonumber && s=s', \omega=\omega '=-\phi, (\tau, \Omega) = (\tau ', \Omega ') \text{ such that } \\
&& (\tau ', \Phi, \Omega ') \notin N^*_{(s, -\omega, \omega)}(\mathbb{B})\}
\end{eqnarray}
We set $L_{\mathbb{B}} := C_{\mathbb{B}} \circ \hat{L}$ and view this as our new ``bad" set. As $\hat{\Lambda}_+ \subset N^*(\{s=0, \phi=\omega\})$ it follows that $C_{\mathbb{B}} \circ \hat{\Lambda}_+ = \emptyset$.  It is easy to verify the compositions with our reflected Lagrangians are again transverse, yielding 
\begin{eqnarray} 
\label{backscatparam} \Lambda^a_{\mathbb{B}} &=& \{(-2y \cdot \frac{\nu}{|\nu|}, \frac{\nu}{|\nu|}; |\nu|, |\nu|i^*_{\frac{\nu}{|\nu|}}(y) : (y, \nu) \in N^*(S_1) \} \\
\Lambda^c_{\mathbb{B}} &=& \{(-2y \cdot \frac{\nu}{|\nu|}, \frac{\nu}{|\nu|}; |\nu|, |\nu|i^*_{\frac{\nu}{|\nu|}}(y) : (y, \nu) \in N^*(S_2) \}.
\end{eqnarray}

A calculation similar to those in Section \ref{lscr} shows $\Lambda^a_{\mathbb{B}} \cap \Lambda^c_{\mathbb{B}}$ is again another codimension $d_2$ submanifold of both $\Lambda^a_{\mathbb{B}}$ and $\Lambda^c_{\mathbb{B}}$.  The polar coordinate argument in the final part of Section \ref{lscr} proves this intersection is clean.  We apply Theorem \ref{fiomapping} from Section \ref{pwo} to get the first part of our main result.

\begin{thm} \label{finalthm} The backscattering, the full scattering kernel restricted to the backscattering data, is
\begin{eqnarray*}
\nonumber \alpha_{| \mathbb{B}} &=& j_{\mathbb{B}}^*(\alpha) \in \\
 &&I^{\frac{n-1}{4} + M_1 + M_2 + \frac{d_1 + d_2 + 1}{2} - \frac{n}{2} - \frac{1}{2} + \frac{2n-1}{4}, -M_2 - \frac{d_2}{2}}(\Lambda^a_{\mathbb{B}} \backslash (L_{\mathbb{B}} \cap \Lambda^a_{\mathbb{B}}), \Lambda^c_{\mathbb{B}} \backslash (L_{\mathbb{B}} \cap \Lambda^c_{\mathbb{B}})) \\
\nonumber &+& I^{\frac{n-1}{4}+ M_1 + \frac{d_1 + 1}{2} - \frac{n}{2} - \frac{3}{2} + \frac{2n-1}{4}}(\Lambda^a_{\mathbb{B}} \backslash (L_{\mathbb{B}} \cap \Lambda^a_{\mathbb{B}})) \\
\nonumber &+& I^{\frac{n-1}{4} + M_1 + M_2 + \frac{d_1 + d_2 + 1}{2} - \frac{n}{2} - \frac{3}{2} + \frac{2n-1}{4}}(\Lambda^c_{\mathbb{B}} \backslash (L_{\mathbb{B}} \cap \Lambda^c_{\mathbb{B}})) \\
&+& H_{L_{\mathbb{B}}}^{s_0 - \frac{2n-1}{4} - \frac{n-1}{4}} + H^{\tilde{s} - \frac{2n-1}{4} - \frac{n-1}{4}}
\end{eqnarray*}
and determines the submanifolds $S_1$ and $S_2$ as well as $\mu(q)$, for $s_0$ and $\tilde{s}$ as in Theorem \ref{comparison}.
\end{thm}

\noindent \underline{Proof}:  Consider the parametrization of $\Lambda^a_{\mathbb{B}}$ in (\ref{backscatparam}). Then $\frac{\Omega}{\tau}-\frac{1}{2}s \omega = i^*_{\frac{\nu}{|\nu|}}(y) + (y\cdot\frac{\nu}{|\nu|})(\frac{\nu}{|\nu|}) = y$.  This means that $S_1$ is determined by $\Lambda^a_{\mathbb{B}}$.  The same holds for $S_2$ from the parametrization of $\Lambda^c_{\mathbb{B}}$.  

Theorem \ref{comparison} tell us the principal symbol $\mu(\alpha_{|\mathbb{B}})$ of the exact scattering kernel is equal to $\mu(\rho^* F \square^{-1}(q \cdot \delta))$.  The ellipticity of $F, \rho^*, \text{ and } j_{\mathbb{B}}^*$ implies $\mu(\alpha_{|\mathbb{B}})$ is a non-zero factor times $\mu(\square^{-1}(q\delta))$.   By results on the symbol calculus of FIOs in \cite{Ho71}, we are able to divide out by these multiples and focus on $\mu(\square^{-1}(q\delta))$.

Microlocally near $\Sigma_1 \cup \Sigma_2$, $\square^{-1}$ acts as an elliptic FIO associated to $C_{\square}$ described in (\ref{flowout}).  Recall Theorem \ref{fiomapping} and the discussion after Remark \ref{intuition} which sets up the oscillatory representation of the paired Lagrangian in (\ref{spaces}).  We see $\mu(\square^{-1}(q\delta))$, in the coordinates given by (\ref{flowout1}) or (\ref{flowout2}), is another non-zero elliptic factor times $\sigma(q\delta)$ on $\Lambda^a_{1+} \cup \Lambda_{2+}^c$ by symbol calculus results in \cite{GuiUhl81}.

As $\mu(\square^{-1}(q\delta))(x - r(\nu -\sigma \omega),x \cdot \omega + r, \omega ; \nu -\sigma \omega, \sigma, -\sigma i^*_{\omega}(x))$ is a function whose variables parametrize the bicharactersitics that foliate $\Lambda^a_{1+}$ and $\Lambda_{2+}^c$,  we can flow our symbol back in the $r$ to $\Sigma_1 \cup \Sigma_2$, obtaining $\mu(\square^{-1}q\delta)=\mu(\square^{-1})(x,x \cdot \omega, \omega ; \nu -\sigma \omega, \sigma, -\sigma i^*_{\omega}(x)) \times \mu(q)(x,\nu) \times \sigma(\delta)(x,x \cdot \omega,\omega;\nu, \tau, \Omega)$ near $\Sigma_1 \cup \Sigma_2$.  Note that $\mu(\delta)(x,t,\omega;\nu, \tau, \Omega)=1$ and $\mu(\square^{-1})(x, x \cdot \omega, \omega ; \nu -\sigma \omega, \sigma, -\sigma i^*_{\omega}(x)) \neq 0$ in this region.  Therefore, we can divide out by these elliptic factors, leaving us with $\mu(q)(x,\nu)$.  As $\Sigma_1 \cup \Sigma_2$ is a dense open subset of $\Sigma^1 \cup \Sigma^2$, it is possible to recover $\mu(q)$ on all of $\Lambda_1 \cup \Lambda_2$ as $\mu(q)$ is function of only $x$ and $\nu$, with no restrictions on either variable besides lying in $\Sigma_1 \cup \Sigma_2$. $\blacksquare$

It is not hard to generalize this result.  In fact, the operator $j^*_{\mathbb{D}}$ that restricts to a submanifold $\mathbb{D}$ is another elliptic Fourier integral operator that satisfies the necessary transversality conditions.  If we have the Lagrangians associated to our new scattering data $\alpha_{| \mathbb{D}}$ cleanly intersecting in codimension $d_2$, Theorem \ref{finalthm} still holds.  Since this last step follows immediately from the general form of Theorem \ref{fiomapping} that appears in \cite{GuiUhl81}, we have
\begin{cor}
If $\mathbb{D} = \{(s,\phi,\omega): \phi = \varphi(s,\omega)\}$ with
\begin{enumerate}
\item
$\varphi(s, \omega) \neq \omega$ for all $\omega \in S^{n-1}$,
\item
$\varphi_s(\omega)=\frac{\omega - \varphi(s, \omega)}{|\omega - \varphi(s, \omega)|}$ is an automorphism of $S^{n-1}$ for all $s \in \mathbb{R}$,
\item
$\varphi^{-1}(s, \nu) \cdot \nu \neq 0$ for all $\nu$ in the image of the Gauss map of $S_1$ and $S_2$, for every $s \in \mathbb{R}$,
\end{enumerate}
then $\alpha_{|\mathbb{D}}$ determines $S_1$ and $S_2$ as well as $\mu(q)$.
\end{cor}
\begin{rem}
Conditions 1-3 are sufficient to carry out a calculation similar to that in the proof of Theorem \ref{finalthm} describing the steps to reconstruct the submanifolds $S_1$ and $S_2$.
\end{rem} 

\section{Appendix: Multiphases}

\noindent The following is a proposition from \cite{Mend82}.
\begin{prop} \label{multphasprop}
Let $\tilde{\Lambda}_0, \tilde{\Lambda}_1$ be Lagrangian submanifolds of $T^*(X) \backslash 0$, $\lambda_0 \in \tilde{\Lambda}_0, \tilde{\Lambda}_1$ and $p_1$ be a homogeneous function of degree 1 such that $p_1(\lambda_0)=0$ and the Hamiltonian vector field $H_{p_1}$ associated to $p_1$ is not tangent to $\tilde{\Lambda}_0$.  If $\tilde{\Lambda}_1$ is the flowout from $\tilde{\Lambda}_1 \cap \{p_1=0\}$ by $H_{p_1}$, then there is a multiphase function $\varphi$ that parametrizes the pair $(\tilde{\Lambda}_0, \tilde{\Lambda}_1)$ which can be chosen such that $\frac{\partial \varphi}{\partial s}(x,s, \theta) = p_1(x, d_x\varphi)$ and $\varphi(x,0,s)=\varphi_0$ with $\varphi_0$ a phase function parametrizing $\tilde{\Lambda}_0$. 
\end{prop}
\begin{rem}
Proposition \ref{multphasprop} can be generalized to the situation where the Lagrangians cleanly intersect in a higher codimension.
\end{rem}
\noindent We will use this proposition to determine a multiphase functions for the pairs $(\Lambda_{1+}, \Lambda_{1+}^a)$ and $(\Lambda_{2+}, \Lambda_{2+}^c)$ that appear in Section \ref{bose}.  Let us compute the phase for first pair $(\Lambda_{1+}, \Lambda_{1+}^a)$.  By the proposition of Mendoza, obtaining the desired multiphase involves solving the following initial value problem,
\begin{equation}
\begin{cases} \label{hamjac}
\frac{\partial \varphi}{\partial s} - p(d_x\varphi, \frac{\partial \varphi}{\partial t}) = 0 \\
\varphi(x,t, \omega; \theta, \sigma, 0) = \varphi_0
\end{cases}
\end{equation}
where $p(\xi, \tau) = \frac{|\xi|^2}{\tau} - \tau$ defines the characteristic variety of $\square$ in $T^*(\mathbb{R}^n \times \mathbb{R} \times S^{n-1})$ and $\varphi_0(x,t,\omega; \theta, \sigma) = \vec{H}(x) \cdot \theta + (t - x \cdot \omega)\sigma$.  We are subsequently led to the system of ordinary differential equations 
\large{\begin{equation}
\begin{cases}
\frac{dx_i}{dr} = -2 \frac{\xi_i}{\tau}, & \text{ for } i=1,...,n \\
\frac{dt}{dr} = -\frac{|\xi|^2}{\tau^2} - 1 \\
\frac{d \xi_i}{dr} = 0, & \text{ for } i=1,...,n \\
\frac{d \tau}{dr} = 0 \\
\frac{ds}{dr} = 1 \\
s(0) = 0
\end{cases}
\end{equation}}
\normalsize
Solving for the characteristics, after replacing $r$ for $s$, we come to 
\begin{eqnarray}
x_i(s)&=&(-2\frac{\xi_i}{\tau})s + x_i^o, \text{ for }i=1,...,n \\
t(s) &=& (- \frac{|\xi|^2}{\tau^2} - 1)s + t^o.
\end{eqnarray}
Now solving for $x_i^o$ and $t^o$ in terms of $s$ and plugging the resulting equations into $\varphi_0$ yields
\begin{equation}\label{compmultphase}
\varphi = \vec{H}_1(x + \frac{2d_x \varphi_0}{\sigma}s)\cdot \theta + (t + (\frac{|d_x \varphi_0|^2}{\sigma^2}+1)s - (x - \frac{2d_x \varphi_0}{\sigma}s \cdot \omega))\sigma
\end{equation}
as our desired multiphase.  Substituting $\vec{H}_2$, the defining functions for $S_2$, for $\vec{H}_1$ in (\ref{hamjac}) and (\ref{compmultphase}) gives us a multiphase for the pair $(\Lambda_{2+}, \Lambda_{2+}^c)$.

\vskip .2cm

\noindent {\bf{Acknowledgements}:}  This article is a revision of the author's Ph.D. thesis.  The author is most grateful to his advisor, Allan Greenleaf, for his help and guidance.  

A portion of this work was supported by NSF grant DMS-0853892.

\bibliographystyle{elsarticle-harv}
\bibliography{references}

\end{document}